\documentclass[preprint,authoryear,12pt]{elsarticle}
\usepackage{amssymb}
\usepackage{color}
\usepackage{amsmath}
\usepackage{subfigure}
\usepackage{graphicx}
\usepackage{graphics}
\usepackage{booktabs}
\usepackage{threeparttable}
\usepackage{appendix}
\usepackage{pdflscape}
\usepackage{epsfig}
\usepackage{booktabs}
\usepackage{bbding}
\usepackage{multirow}
\usepackage{geometry}
\usepackage[displaymath]{lineno}
\usepackage{setspace}
\usepackage[colorlinks,linkcolor=red,anchorcolor=blue,citecolor=green]{hyperref}
\usepackage{soul}
\usepackage{subfigure} 
\usepackage{longtable}
\usepackage{mathpazo}
\usepackage{appendix}

\makeatletter

\newcommand{\Rmnum}[1]{\expandafter\@slowromancap\romannumeral #1@}
\makeatother

\journal{Transportation Research Part A}

\begin{document}

\begin{frontmatter}

\title{Assessing the dynamic vulnerability of an urban rail transit system and a case study of Beijing, China}

\author[rvt1]{Shouzheng Pan}

\author[rvt1]{Jia He}

\author[rvt2]{Ning Jia}

\author[rvt3]{Der-Horng Lee}

\author[rvt1]{Zhengbing He\corref{cor1}}
\ead{he.zb@hotmail.com}

\address[rvt1]{Beijing Key Laboratory of Traffic Engineering, Beijing University of Technology, Beijing, China}
\address[rvt2]{Institute of Systems Engineering, Tianjin University, Tianjin, China}
\address[rvt3]{ZJU-UIUC Institute, Zhejiang University, Zhejiang, China}

\cortext[cor1]{Corresponding author}

\begin{abstract} 
\begin{spacing}{1.3}
\small
Urban rail transit (URT) is the backbone of a city.
It is significant to understand its vulnerability, i.e., the variation in capacity and service levels of resilience when facing operational disturbances.
Although the network topology is generally fixed, the hourly-changing travel demand greatly impacts on the actual performance of a URT system. 
Unfortunately, few existing studies take such dynamic travel demand into account and ignore its impact on the vulnerability of a URT system. 
To fill the gap, this paper proposes a network vulnerability assessment method with the joint consideration of static network topology and dynamic travel demand, which includes an accessibility-based identification of station importance with time-varying passenger demand and a new dynamic vulnerability evaluation index.
An empirical analysis was carried out by taking the URT system of Beijing, China as an example, and the impact of the more realistic multiple consecutive station failures in a URT system is also examined.
Results show that the distribution of high-importance stations indeed varies with the time of day, affected by both static topology and hourly-changing passenger flow. 
When the disturbance of operation delay occurs, the impact of high-importance stations on the network vulnerability changes nonlinearly with the increase of delayed travel demand. 
Some stations that serve as bridges and are visited by large passenger flows have the greatest impact on network vulnerability. 
Network performance degradation is obviously segmented and stratified in the case of interval continuous failure. 
The disruption between different lines is the main cause of network performance degradation, and some high-importance stations within the lines act as catalysts to accelerate the performance degradation. 
The proposed method provides a reference for measuring dynamic passenger flow-related network vulnerability and supplies the field with a new vulnerability evaluation index.

\end{spacing}
\end{abstract}


\end{frontmatter}


\begin{spacing}{1.4}

\section{Introduction}

Urban rail transit (URT) is one of the most important transportation systems in a city, covering almost every corner in a city such as New York City, Paris, London and Beijing. 
The safe and reliable operation of a URT system is a guarantee for people's traveling and social activities. 
However, a URT system is occasionally impacted by external and internal disturbances, such as extreme weather, large events, power and signal failures.
The disturbances may lead to train delays, station closure, or even the suspension of a whole line. 
It is significant to understand the performance changes in a URT system under disturbances, including the identification of key components, network vulnerability, magnitude of the impact of the disturbance, etc.

As an important component of resilience, vulnerability quantifies the speed and extent of performance degradation \citep{Reggiani2015a}. 
It is one of the most relevant concepts for describing transportation system performance under disturbances, i.e., \textit{susceptibility to unusual incidents that can result in considerable reductions in system serviceability} \citep{Berdica2002}.
Regarding the vulnerability of URT, early studies focused mainly on static network topology.
For example, \cite{Crucitti2003} analyzed the characteristics of URT topology, including node degree, betweenness, and global efficiency, and found that most URT networks show small-world characteristics and have better connectivity reliability. 
However, the distribution of travel demand in a city is well known to vary with the time of day, resulting in a time-varying imbalance between transit supply and travel demand.
Therefore, the performance of a URT system is dynamic, making the na\"ive analysis of static network topology insufficient to assess a URT system with time-varying performance.

Until recently, passenger flow was introduced into the studies of the vulnerability of a transportation network, such as the rail transit \citep{Yang2015, Yin2016}, road transportation \citep{Chen2012}, and public transportation \citep{Rodriguez-Nunez2014,Cats2014}.  
For example, \cite{Feng2017} proposed a two-layer network model to analyze the vulnerability of a URT network under the influence of passenger flow. 
The results show that the flow network and train network have similar spatiotemporal clustering characteristics and are both affected by network topology. 
To reflect the impact of passenger flow on system performance, various vulnerability metrics were established by constructing a topological and weighted network model combining complex network theory and passenger flow data \citep{Zhou2019}. 
\cite{Ouyang2015} analyzed the vulnerability of railway and airline systems using daily accessibility as metrics, which represents the proportion of nodes that are reachable to each other in a network.
\cite{Hong2020} used daily changes in accessibility indicators to calculate the vulnerability of integrated URT and high-speed rail systems.
The distribution of high-importance stations under different evaluation indices was found to be different.
After that, \cite{adjetey2016} constructed a simulation model to measure the resilience in a railway system combining the changes of passengers flow and delayed volume as metrics.
\cite{Voltes-Dorta2017} also adopted the total passengers delay in one day as metrics. 
The results show that the metrics considering passenger recovery and demand changes at different times are more accurate than topological characteristics in analyzing transportation system performance.
\cite{Lu2018} established a resilience evaluation index for daily operational incidents by considering both delayed origin-destination (OD) volume and network topology. 
Except for time difference, \cite{Ganin2019} evaluated travel demand for delay in different regions and scenarios on the intelligent transportation system.
To capture the impact of dynamic demands, a combination of metrics is also often used.
\cite{Knoop2012} proposed a comprehensive metrics and considered the total travel time by conducting linear fitting of multiple vulnerability indices.
However, the metrics or methods proposed in most of the existing studies are still relatively static, i.e., network weights are measured using passenger flow that was collected in a relatively large time period (e.g., a whole day).
The granularity is unable to accurately reflect time-varying travel demand within a day (e.g., peak-and-off-peak patterns) as well as its impact on network vulnerability.
One of the obstacles might be the high computational efficiency required by frequently updating the route and line passenger flow every short time interval.

In addition, different disturbances usually have different impacts on network vulnerability.
It is necessary to distinguish disruption types and analyze system performance in specific disturbance scenarios \citep{Donovan2017, Ilbeigi2019, Diab2019}.
Different transportation systems often face different kinds of disturbances, such as the traffic congestion on a road network \citep{khaghani2019,Zhang2019b} and flight delay in an aviation network \citep{Voltes-Dorta2017}.
For a URT system, train delays that may result from, e.g., natural disasters and operation incidents, aggravate the imbalance between demand and supply and degrade system performance.  
\cite{Lu2018} analyzed the impact of train delays on system performance under operational incidents, showing that network resilience gradually decreases as train delay increases and the robustness of transfer stations is better than that of nontransfer stations.
The failure of stations or lines is usually taken as a typical trigger of operational delay \citep{Zhang2019}.
For example, \cite{Sun2016} analyzed the vulnerability of a URT network from the perspective of operational lines, instead of taking traditional nodes or connecting edges as interference targets.
The results show that passenger flow is the key factor that causes network vulnerability, and the failure of loop lines has the greatest impact on the performance of the whole network.
However, the failure of a whole line might rarely occur in a modern URT system with advanced operation and management technologies. 
A more practical scenario is the failure of multiple consecutive stations at the same time, which is also known as {\it interval failure}. 
For example, a train rear-end incident occurred in Shanghai Metro, China due to equipment failure on September 27, 2011. 
Only the nearest 9 stations along the entire line were closed.
To the best of our knowledge, few studies have focused on the impact of interval failure in a URT system.

In summary, although travel demand obviously changes within a day and may have a significant impact on network vulnerability during different time periods (e.g., morning and noon), the existing literature pays less attention to it.
Moreover, the impact of the more realistic multiple consecutive station failures in a URT system is still unclear. 
To fill these gaps, this paper proposes a dynamic vulnerability assessment method to consider the impact of short-delay and long-delay disruptions on system performance under time-varying travel demand. 
The station and interval failures are explicitly tested by taking the URT system in Beijing, China as a case study.
The main contributions of this paper are as follows.
\begin{itemize}
	\setlength{\itemsep}{0pt}
	\setlength{\parsep}{0pt}
	\setlength{\parskip}{0pt}	
  \item An accessibility-based identification method for station importance is proposed by jointly considering static network topology and dynamic travel demand. 
  Instead of dynamic passenger flow assignment, it is calculated based on merely passenger demand that has direct connections with the disrupted stations, making it possible to assess network vulnerability in a real-time manner. 
  \item A new dynamic vulnerability evaluation index is proposed to assess network vulnerability changes under short-delay disruptions that result from small disturbances such as train delay and signal failure.
  \item The variation in network vulnerability under the scenarios of station and interval failures is explicitly studied by taking the impact of operation delay into account. To the best of our knowledge, this is the first work that evaluates a URT system under the impact of interval failure.
\end{itemize}

The reminder of the paper is organized as follows.
Section \ref{sec:Methodology} introduces basic definitions and indices that will be used to describe the characteristics of a URT network.
Section \ref{section3} proposes the accessibility-based identification method for station importance and the new dynamic vulnerability evaluation index under disruptions.
By taking the Beijing URT system in China as an example, Section \ref{casestudy} conducts a dynamic vulnerability analysis under the scenarios of station, interval and line failures. 
Finally, summary and discussion of future research are presented in Section \ref{section5}.


\section{Basic definitions of URT network characteristics}\label{sec:Methodology}

\subsection{URT network modeling}\label{sec:network}

A URT network can be abstracted as an undirected and weighted graph composed of nodes and connecting edges.
Specifically, let $\textbf{N}=\{ 1,2,...,N\}$ be a set of nodes (i.e., URT stations composed of transfer and nontransfer stations) with size $N$ and $i\in \textbf{N}$ be a node.
Let $\textbf{E}=\left\{e_{ij}|e_{ij}=\left \langle i, j \right \rangle, i\neq j \right\}$ be a set of edges between nodes $i$ and $j$.
When station $i$ directly connects with station $j$, $e_{ij}=1$; otherwise, $e_{ij}=0$.
Let $\textbf{W}=\left\{w_{ij}|w_{ij}\geq 0\right\}$ be the weight of edge $e_{ij}$, which can be the actual running time $\tau_{ij}$ of trains between stations $i$ and $j$ or the passenger flow $f_{ij}$ at a given time.
Correspondingly, $\textbf{G}=\left\{\textbf{N},\textbf{E},\textbf{W}\right\}$ is the graph.

\subsection{Flow changes with dynamic demand}
Strictly speaking, travel demand in a URT system unevenly changes with time, resulting in a time-varying spatiotemporal distribution of passenger flow inside a URT system.
Loading with such dynamic flow, a URT system that is static in topology shows time-varying characteristics.
To accurately capture the characteristics, it is necessary to assign passenger flow into routes and links, meaning that passenger choice of travel route inside a URT system must be involved.
Unfortunately, the smart card data collected by the most automatic fare collection (AFC) systems cannot reflect the detailed passenger paths inside a URT system.
To remedy this issue, we define a {\it reasonable path} by jointly considering travel time, transfer time, and the number of transfer stations.
The detailed steps are as follows.
(i) Take the train running time as the travel time of passengers between stations.
(ii) Construct a URT topology with virtual transfer arcs according to the actual URT network.
Transfer time is assigned to the virtual arc as its weight.
The number and the weight sum of the transfer arcs along a trip are the transfer number and the total transfer time that a passenger experiences, respectively.
(iii) Calculate the $K$ shortest path between any nodes by using the Dijkstra algorithm.
(iv) Finally, take the path with the minimum total travel time, the smallest transfer number and the minimum transfer time as a reasonable path.
Note that if there are multiple reasonable paths between an OD pair, the paths will be assigned in proportion to their respective travel time.

\textbf{Definition 1}. Among all the shortest paths of any OD pair in a fully connected URT network, the path with the shortest travel time, the smallest transfer number and the shortest total transfer time is a {\it reasonable path} and there is at least one path.

The OD matrices within a time period, which are the reflection of dynamic flow distributions, can be further obtained after possessing a reasonable path between any stations.
The characteristics of a URT system will be analyzed on this basis.

\subsection{Topological characteristics of a URT network with dynamic demand}
Understanding the topological characteristics of a URT network is beneficial to identifying high-importance stations and analyzing network vulnerability. 
This paper first introduces the following three indices that were proposed based on the complex network theory to describe the topological characteristics of a URT network with the consideration of dynamic passenger flow.

\begin{itemize}
	\setlength{\itemsep}{0pt}
	\setlength{\parsep}{0pt}
	\setlength{\parskip}{0pt}	

\item {\it Node degree weighted by passenger flow} (also known as intensity) calculates the passenger flow distribution of a network through passenger flow assignment to reflect the local importance of a network under the influence of travel demand. 
Let $D_i(t)$ denote the degree of station $i$ during time period $t$, and the equation can be written as follows.
\begin{equation}\label{eq1}
  D_i(t)=\sum_{j \in \textbf{N}_{i}} w_{ij}(t)  
\end{equation}
where $w_{ij}(t)$ is the weighted flow from node $i$ to node $j$ during time period $t$ and $\textbf{N}_{i}$ is the set of the nodes that directly connect with node $i$.

\vspace{3mm}

\item  {\it Flow betweenness} is defined as the ratio of the passenger flow passing through station $i$ along a reasonable path to the total passenger flow at a given time period $t$, reflecting the importance of a URT station in terms of being a bridge in the whole network.
Let $B_i(t)$ denote the betweenness during time period $t$, and the equation can be written as follows.
\begin{equation}\label{eq2}
  B_i(t)=\sum_{j=1}^{N} \sum_{k=1}^{N} \frac{f^{i}_{jk}(t)}{f_{jk}(t)},\ \ \ j\in \textbf{N}, \ k\in \textbf{N}, \ i\neq j\neq k 
\end{equation}
where $f_{jk}(t)$ is the passenger flow between any two nodes $j$ and $k$ on a reasonable path  during time period $t$; 
$f^i_{jk}(t)$ is the passenger flow passing through node $i$ along a reasonable path of two nodes $j$ and $k$ during time period $t$.

\vspace{3mm}

\item {\it Demand closeness centrality}, which is defined as the reciprocal of the ratio of the inflow of station $i$ to the total number of stations during time period $t$, reflects the closeness of the distance and travel demand between stations $i$ and $j$. 
Let $C_i(t)$ denote the centrality of station $i$ and the equation can be written as follows.
\begin{equation}\label{eq3}
  C_i(t)=\frac{N-1}{\sum_{j=1}^{N_i(t)} f_{ij}(t)}   
\end{equation}
where $N_i(t)$ is the number of stations where passengers who take a URT train from station $i$ during time period $t$ will arrive in the short future of a URT train leaving station $i$;  
$f_{ij}(t)$ is the passenger flow between stations $i$ and $j$; 
the larger $C_i(t)$ is, the smaller the closeness of station $i$ is.

\end{itemize}

\section{Dynamic assessment of the vulnerability of a URT network } \label{section3}

\subsection{Accessibility-based dynamic identification of station importance}\label{section3.1}

Stations are one of the most basic and critical components for the vulnerability of a URT network.
The importance of a station is usually affected by the topological structure of a network and the passenger flow in the network.
The existing indices that evaluate station importance from the perspective of topological structure include node degree, betweenness centrality and closeness centrality.
To incorporate the impact of travel demand, passenger flow is usually employed to weight the degree (Equation \ref{eq1}) and betweenness (Equation \ref{eq2}).
As mentioned above, most of the existing studies only consider daily passenger flow, obviously ignoring the rapid changes in passenger flow within a day.
One of the reasons is that assigning fine-granularity passenger flow to a network (in particular to a large-scale network) usually creates a high requirement for computational efficiency, which is a challenging task.

To overcome the difficulty, this paper proposes an accessibility-based station importance identification method that does not require passenger flow assignment.
In detail, suppose that the operation of station $i$ is disturbed due to any unexpected reason during time period $t$, and we have the following accessibility-related variables.
\begin{itemize}
	\setlength{\itemsep}{0pt}
	\setlength{\parsep}{0pt}
	\setlength{\parskip}{0pt}	
  \item The following three groups of passengers that have direct connections with station $i$ will be affected during time period $t$, namely, the passengers whose destination is station $i$ (denoted by $S^i_\text{out}(t)$), the passengers whose origin is station $i$ (denoted by $S^i_\text{in}(t)$), and the passengers who have to pass station $i$ (denoted by $S^i_\text{pass}(t)$). 
  The links of the stations that contain station $i$ as an intermediate station are predetermined using \textbf{Definition 1}, and there is no need to make real-time updates. 
  Therefore, $S^i_\text{pass}(t)$ is the number of passengers who enter and leave the URT system from the two ends of the links of stations during time period $t$.
 
  \vspace{3mm}
  
  \item Accordingly, three groups of stations will be affected, namely, the stations used by $S^i_\text{out}(t)$, the stations used by $S^i_\text{in}(t)$, and the stations used by $S^i_\text{pass}(t)$. The numbers of the three groups of stations are denoted by $N^i_\text{out}(t)$, $N^i_\text{in}(t)$ and $N^i_\text{pass}(t)$, respectively.
  Note that not all stations will be used by passengers during time period $t$, particularly when $t$ is quite short. 
  
\end{itemize}

Generally speaking, the influence of a station on a network, i.e., the importance of a station for a network, is positively correlated with $S^i_\text{out}(t)$, $S^i_\text{in}(t)$, $S^i_\text{pass}(t)$ and $N^i_\text{out}(t)$, $N^i_\text{in}(t)$, $N^i_\text{pass}(t)$. Therefore, we propose the following index to assess the importance of a station for a network, which not only considers passengers who would access the disturbed station $i$ during time period $t$ but also has no need for time-consuming passenger flow assignment.
\begin{equation}\label{eq4}
  \zeta_i(t)=\frac{
  S^i_\text{out}(t)N^i_\text{out}(t) + 
  S^i_\text{in}(t)N^i_\text{in}(t)   +
  S^i_\text{pass}(t)N^i_\text{pass}(t)}
  {S(t) N(t)}   
\end{equation}
where the denominator normalizes the impact using the total passengers (denoted by $S(t)$) that enter the URT system during time period $t$ and the number (denoted by $N(t)$) of stations associated with the passengers.
The values of the input variables in Equation \ref{eq4} can be directly and simply obtained from the smart card data collected by the current AFC system, guaranteeing high computational efficiency.

\subsection{Vulnerability measurement of a URT network}

In terms of the resulting delay, the disruptions associated with a URT system can be classified into two types:  
{\it short-delay disruptions} that result from small disturbances such as train delay and signal failure and {\it long-delay disruptions} that result from large disturbances such as fire disaster, heavy maintenance and flood.
Passengers would choose to wait and keep using the delayed route when a short-term disruption occurs. 
In contrast, passengers may choose another travel mode for a replacement of URT when a long-term disruption occurs, since the route affected by a long-term disruption may not be recovered in a short time.
Assume a disruption will trigger a constant delay (denoted by $\tau_\text{disr}$), and we simply use a threshold (denoted by $\tau_\text{disr}^*$) of the delay to distinguish the short- and long-delay disruptions, i.e., 
a short-delay disruption occurs when the delay $\tau_\text{disr} \leq \tau_\text{disr}^*$, and a long-term disruption occurs when $\tau_\text{disr} > \tau_\text{disr}^*$.

To jointly incorporate the disruption delay and dynamic passenger flow, we propose a new metric $\psi_\text{short} (\textbf{G}, t,\tau_\text{disr})$ to describe the vulnerability of a URT network under short-delay disruptions, i.e., 
\begin{equation}\label{eq5}
  \psi_\text{short} (\textbf{G}, t,\tau_\text{disr})=\frac{\sum_{i=1}^{N} \sum_{j=1}^{N} \tau_{ij} (t,\tau_\text{disr}) f_{ij} (t)}{N(N-1)},\ i\neq j, \ i\in \textbf{N} ,\ j \in \textbf{N} ,\ 0< \tau_\text{disr}\leq \tau_\text{disr}^*   
\end{equation}
where $\tau_{ij} (t)$ is the actual travel time between stations $i$ and $j$ during time period $t$ and is calculated as 
\begin{equation}
    \tau_{ij} (t, \tau_\text{disr}) = \tau_{ij} + \theta_{ij}(t)\tau_\text{disr}
\end{equation}
where $\tau_{ij}$ is the constant travel time between stations $i$ and $j$ under normal conditions;
$\theta_{ij}(t)$ is a dummy variable, i.e., $\theta_{ij}(t)=1$, if the disruption with delay $\tau_\text{disr}$ occurs during time period $t$ and is associated with stations $i$ or $j$; $\theta_{ij}(t)=0$, otherwise.

To measure the network vulnerability under long-delay disruptions, we first remove the disturbance-influenced station and edges and reconstruct the network graph $\textbf{G}'=\{\textbf{N}', \textbf{E}', \textbf{W}'\}$. 
Then, after each disturbance-disruption, we measure vulnerability by employing the following index, which represents the variation in URT network performance. 
\begin{equation}\label{eq6}
  \psi_\text{long} (\textbf{G}', t,\tau_\text{disr})=\frac{\sum_{i=1}^{N'} \sum_{j=1}^{N'}  \frac{1}{\tau_{ij}'}f_{ij}'(t)}{N'(N'-1)},\ i\neq j, \ i\in \textbf{N}' ,\ j \in \textbf{N}' ,\ \tau_\text{disr}> \tau_\text{disr}^*  
\end{equation}
where $N'$ is the number of nodes in the reconstructed network $\textbf{N}'$;
$\tau_{ij}'$ and $f_{ij}'(t)$ are the constant travel time and the passenger flow during time period $\tau$ between stations $i$ and $j$ on a reasonable path of the reconstructed graph $\textbf{G}'$;
$\psi_\text{long} (\textbf{G}', t,\tau_\text{disr})$ is also known as operational efficiency \citep{Sun2016}.


\section{A case study of Beijing URT, China}\label{casestudy}

\subsection{Beijing URT and data description}

By 2017, the URT system in Beijing, China contained 19 lines (Lines 2 and 10 were circle lines), and the total length of the lines had reached 574 kilometers.
The average number of daily passengers was approximately 12.4 million.
It is representative and significant to assess the dynamic vulnerability of such a URT system using the proposed metrics, considering the fact that the Beijing URT is one of the largest URT systems worldwide.

\begin{figure}[!htbp]
    \centering
    \subfigure[Time-varying demands]{
    \includegraphics[width=2.05in]{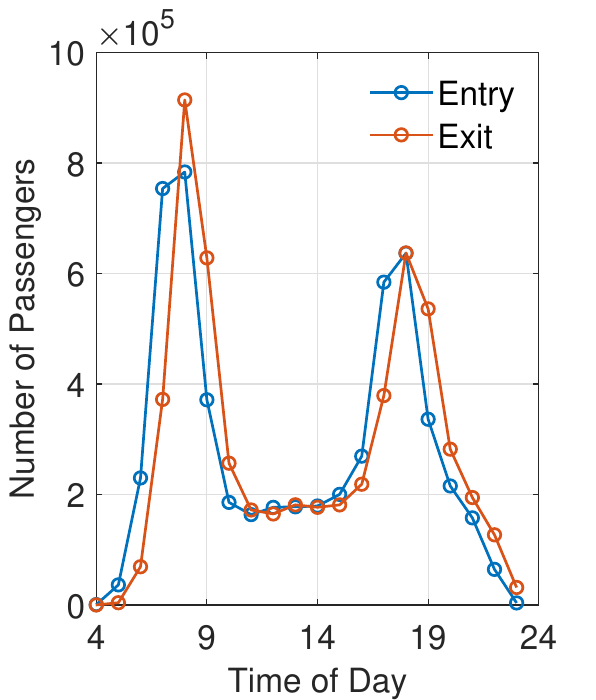}
    \label{fig:f1}
    }%
    \subfigure[Passenger demands for operational lines at different time periods]{
    \includegraphics[width=4.7in]{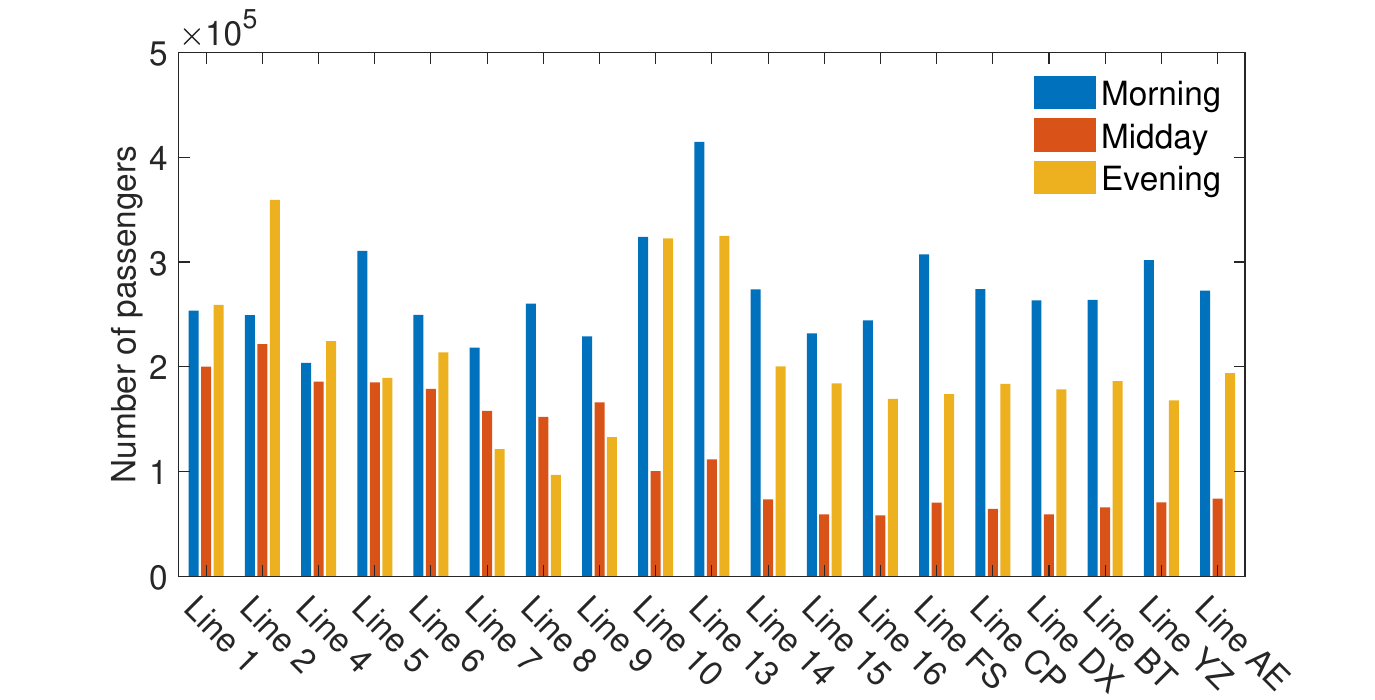}
    \label{fig:f3}
    }
    \caption{Dynamic passenger demands of the Beijing URT system on October 16th, 2017.
    Morning: 6:00 - 8:00; 
    Midday: 11:00 - 13:00; 
    Evening: 16:00 - 18:00.
     FS: FANGSHAN; CP: CHANGPING; DX: DAXING; BT: BATONG; YZ: YIZHUANG; AE: Airport Express
    }
    \label{fig:passenger demand}
\end{figure}

The Beijing URT system is abstracted as an undirected and weighted topological network with 287 stations, including 54 transfer stations, as the way introduced in Section \ref{sec:network}. 
The data used here were collected by the AFC system of the Beijing URT in a typical weekday, i.e., October 16th (Monday), 2017.
After data cleaning, the number of valid records is approximately 5.52 million.
The data contain every passenger's record of entering and leaving the URT system, including the time and the entry and exit stations.
As shown in Figure \ref{fig:passenger demand}, passenger demands of the Beijing URT system at different time periods and operational lines can be calculated.
Figure \ref{fig:f1} presents the hourly change in the number of passengers entering and leaving the URT system.
Peak-and-off-peak patterns in travel demand can be clearly observed.
Figure \ref{fig:f3} presents heterogeneous demands for URT lines, and we have the following observations: (1) different lines correspond to quite different passenger demands, and (2) passenger demands are also different at different time periods for the same line.
The heterogeneous observations of spatiotemporal demand clearly indicate the necessity of the dynamic assessment of the time-changing vulnerability of a URT system under the occurrence of the failure of stations, intervals and even whole lines. 

\subsection{Identification of dynamic importance of stations}\label{section4.2}

\begin{figure}[!ht]
\centering
\subfigure[Type-\Rmnum{1}]{
\label{fig:f4.sub1}
\includegraphics[width=3in]{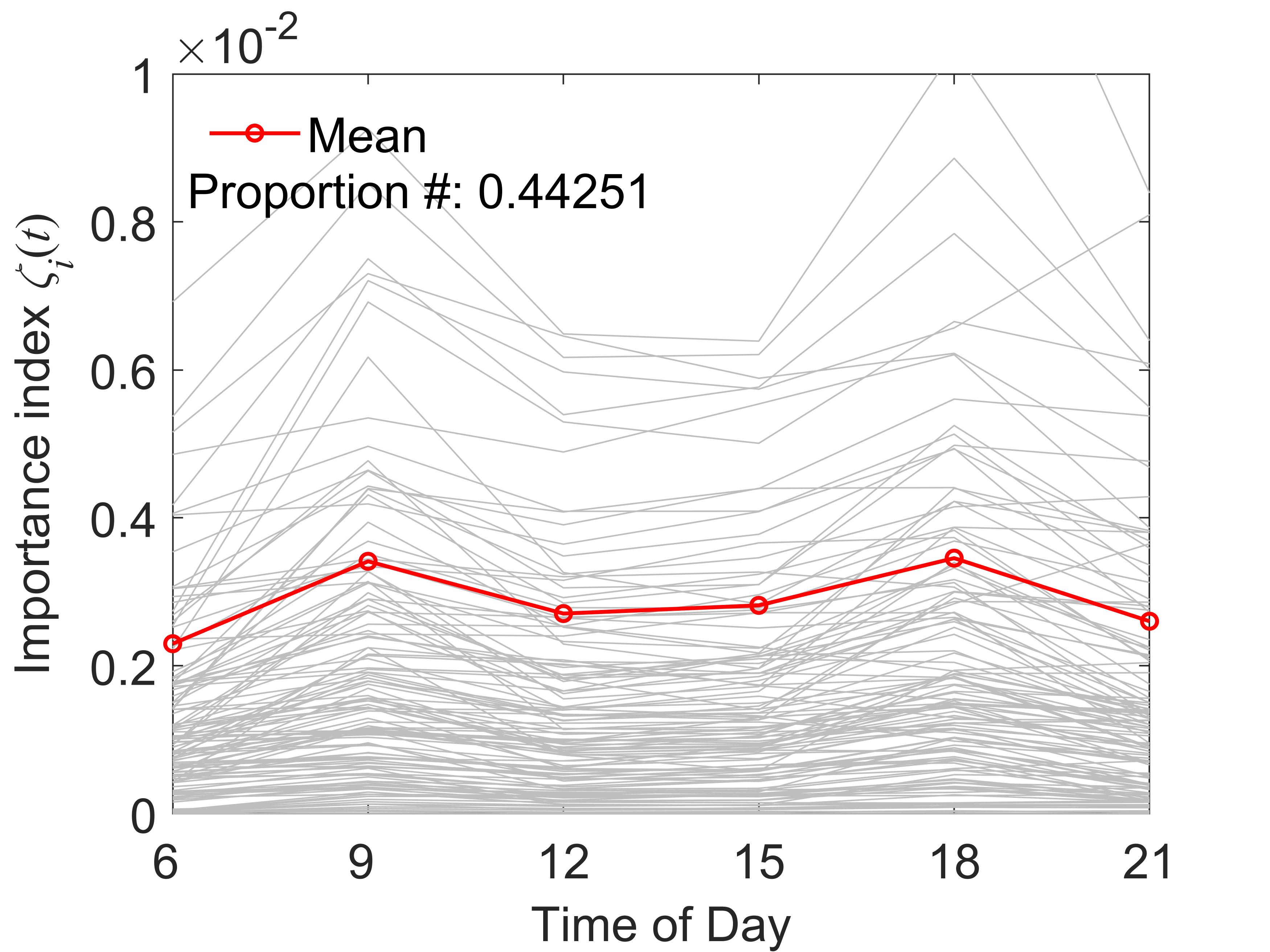}
}
\subfigure[Type-\Rmnum{2}]{
\label{fig:f4.sub2}
\includegraphics[width=3in]{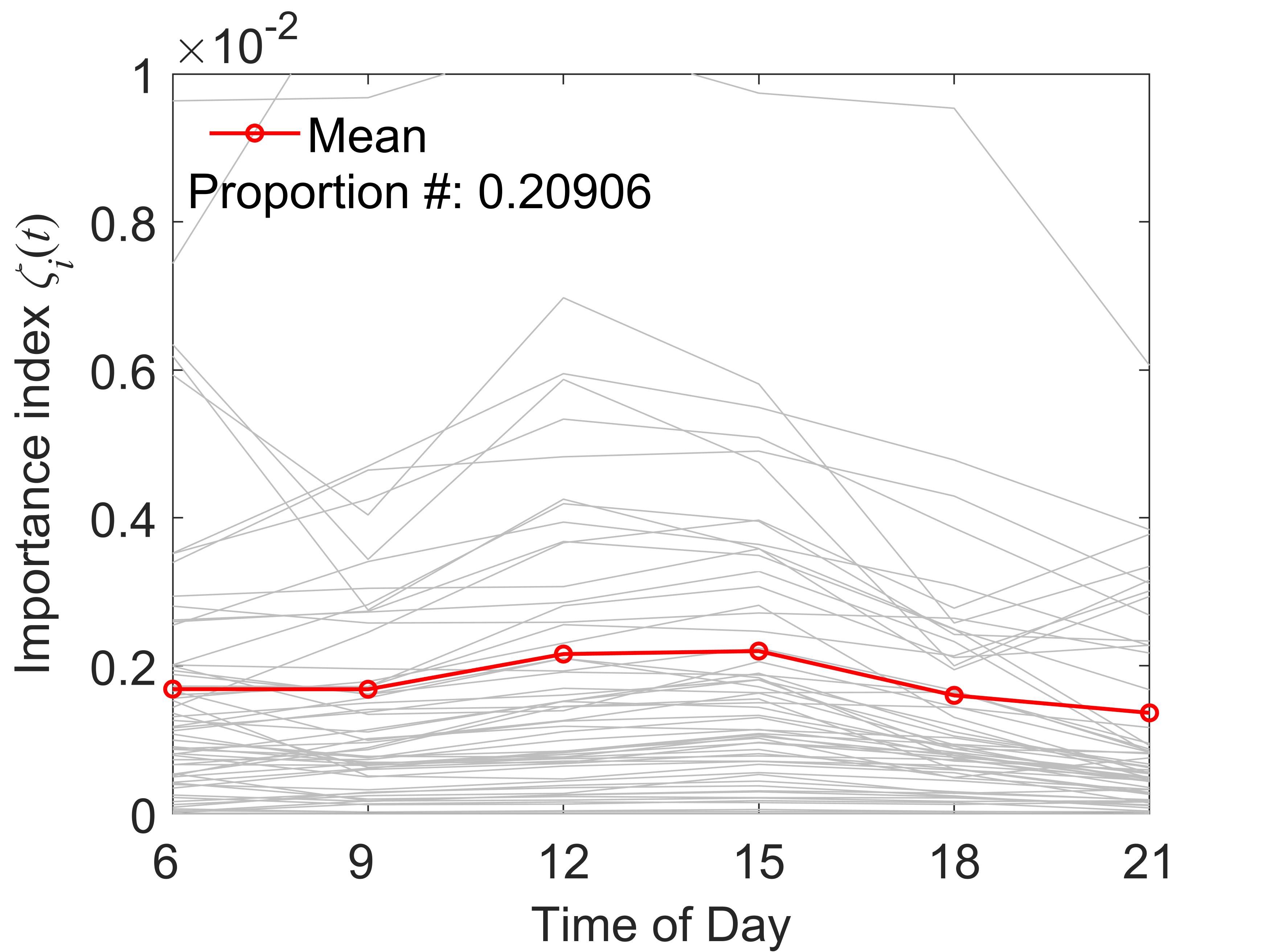}
}

\subfigure[Type-\Rmnum{3}]{
\label{fig:f4.sub3}
\includegraphics[width=3in]{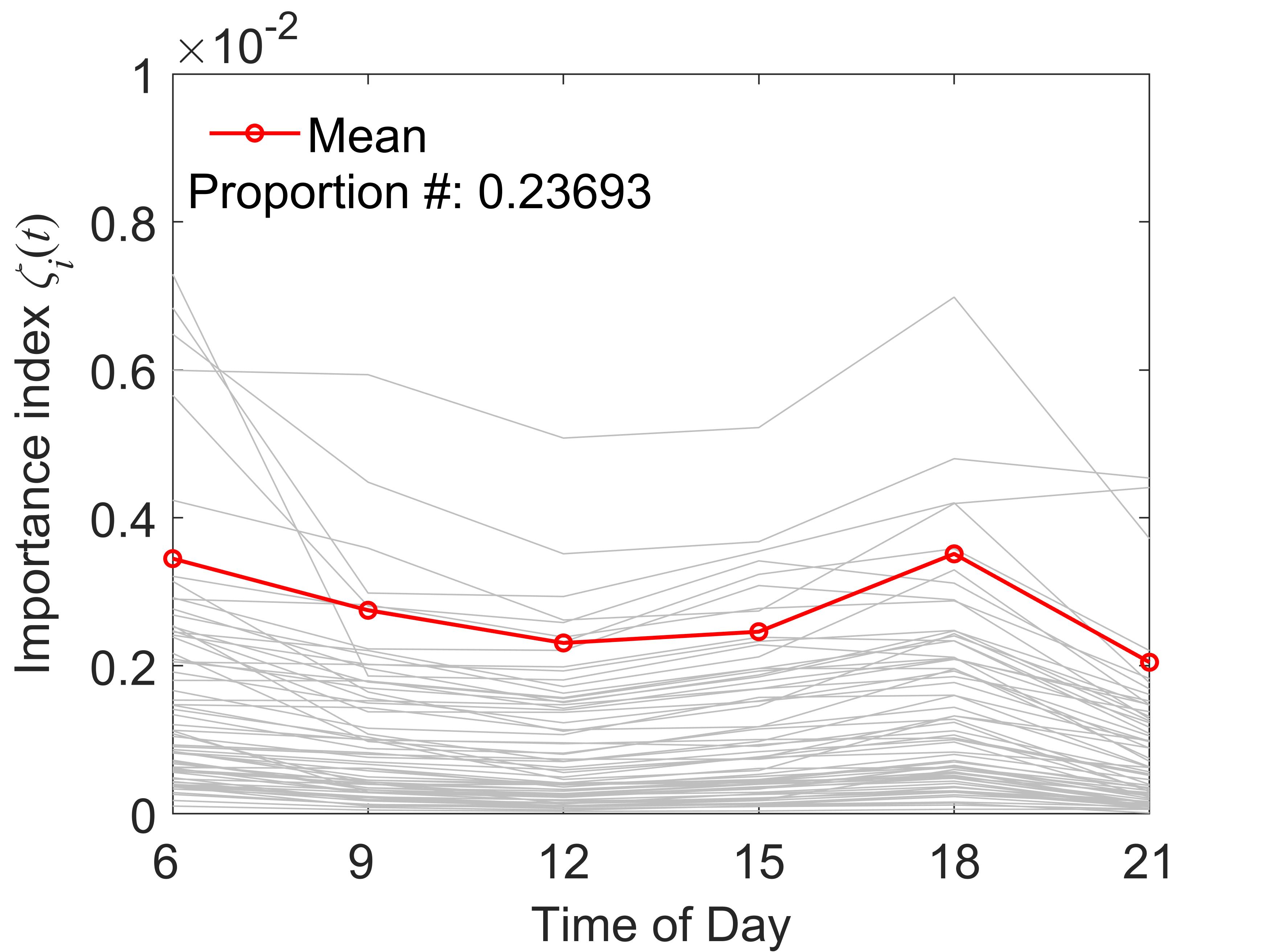}
}
\subfigure[Type-\Rmnum{4}]{
\label{fig:f4.sub4}
\includegraphics[width=3in]{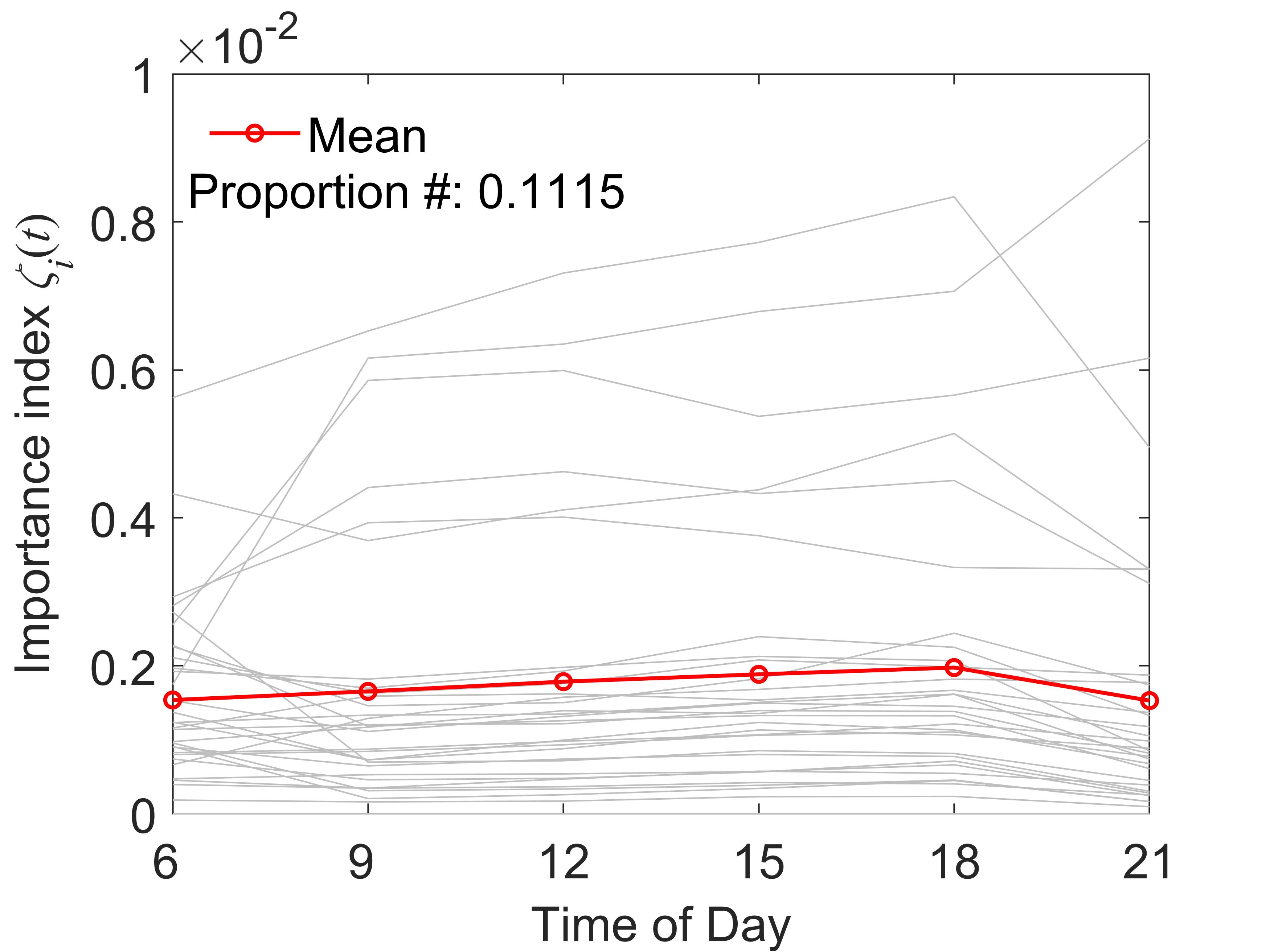}
}
\centering
\caption{Four typical types of changes in station importance. "Proportion" is the ratio of the number of stations of one type to the total number of stations }
\label{fig:f4}
\end{figure}

The dynamic importance of stations of the Beijing URT is identified by employing the proposed accessibility-based index (Equation \ref{eq4}).
The importance of each station at different times of the day can be calculated, generating a continuous curve that changes over time.
There are 287 unique stations in the Beijing URT network, and thus, we obtain 287 curves to represent the time-varying importance of stations.
To understand the changing patterns of these curves, we conduct a cluster analysis on the curves by using hierarchical clustering.
In detail, the importance of each station is calculated every three hours, and the time range is from 5:00 to 23:00, i.e., 18 hours in total. 
The slope between two successive points in the curve is used as the distance for hierarchical clustering, and four types of those curves can be obtained, as shown in Figure \ref{fig:f4}. 
The variation trend and key features are summarized as follows.
\begin{itemize}
	\setlength{\itemsep}{0pt}
	\setlength{\parsep}{0pt}
	\setlength{\parskip}{0pt}
	\item Type-\Rmnum{1}: There are two peaks of importance for Type-\Rmnum{1} stations. High importances occur mainly in the morning and evening periods, i.e., the peak hours of travel demand (refer to Figure \ref{fig:f1}). 
	
	\item Type-\Rmnum{2}: The importance of stations in Type-\Rmnum{2} have a single peak in the midday period. 
	
	\item Type-\Rmnum{3}: There are also two peaks of importance for Type-\Rmnum{3} stations, while the importance peak in the morning period is earlier compared with Type-\Rmnum{1} stations. 
	
	\item Type-\Rmnum{4}: Most of these stations have small importance indices and do not have distinct peaks, which may be related to the topological structure of the station and constant travel demand.
	
\end{itemize}

\begin{figure}[!htbp]
\centering

\subfigure[Type-\Rmnum{1}]{
\label{fig:f5.sub1}
\includegraphics[width=3.3in]{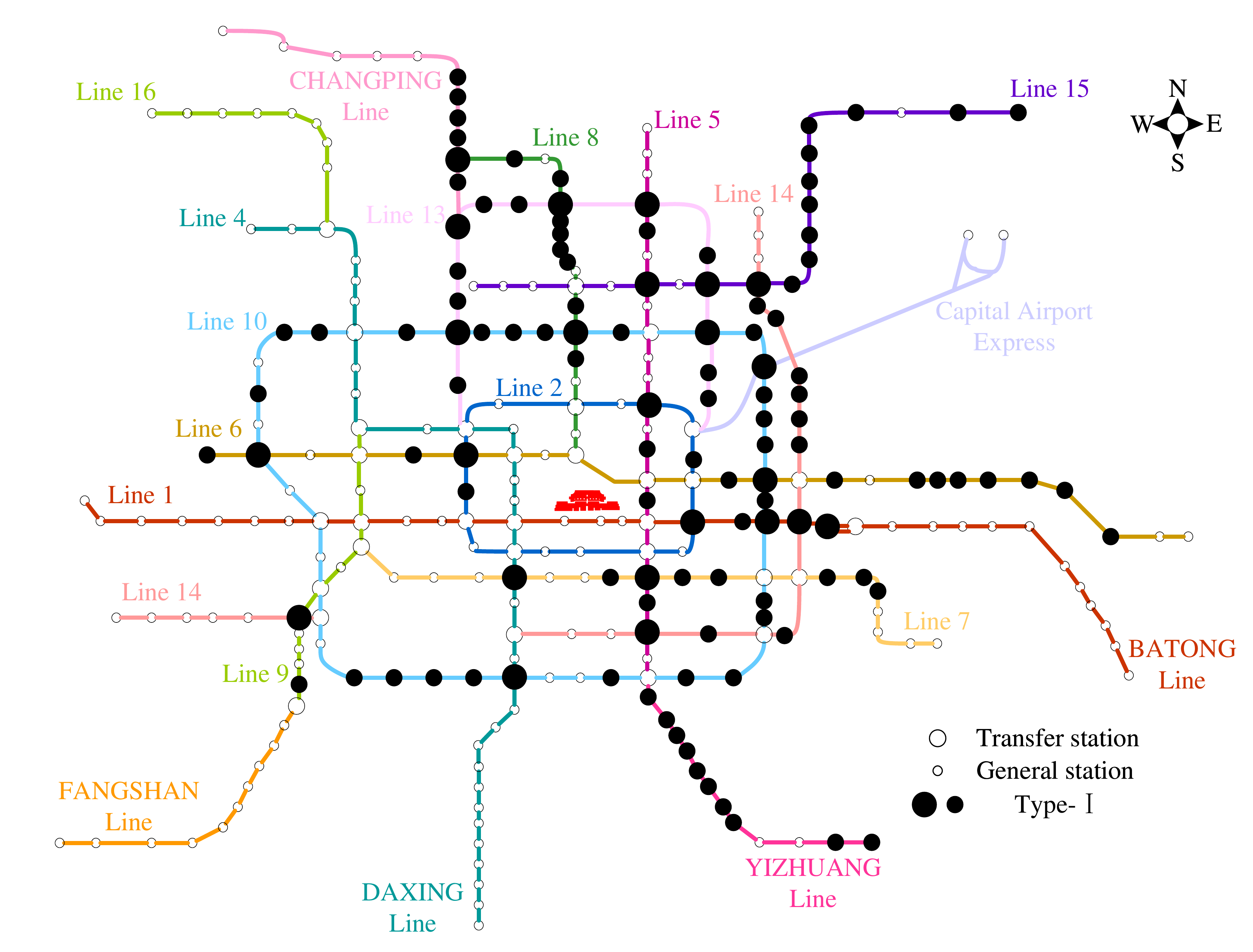}
}%
\subfigure[Type-\Rmnum{2}]{
\label{fig:f5.sub2}
\includegraphics[width=3.3in]{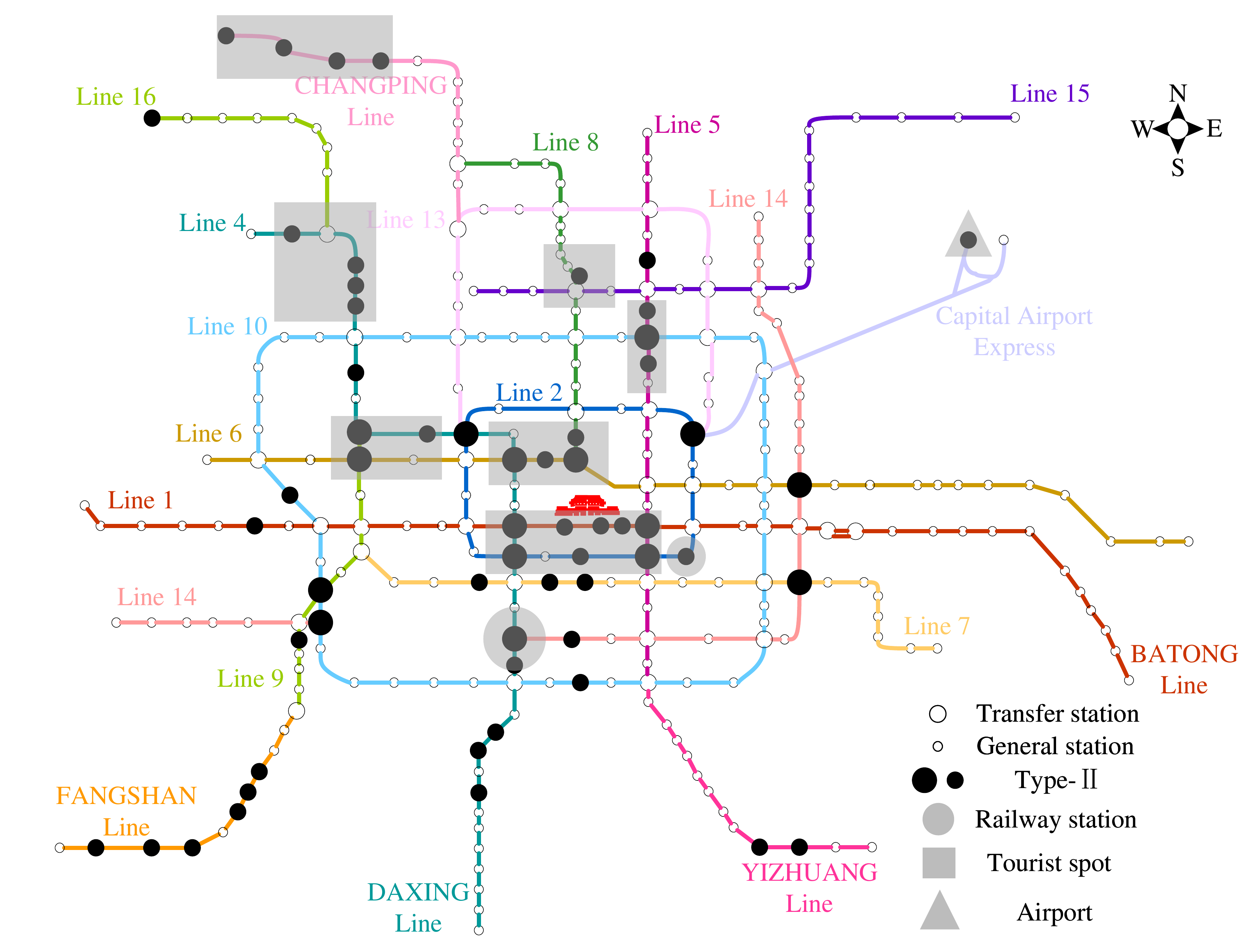}
}%

\subfigure[Type-\Rmnum{3}]{
\label{fig:f5.sub3}
\includegraphics[width=3.3in]{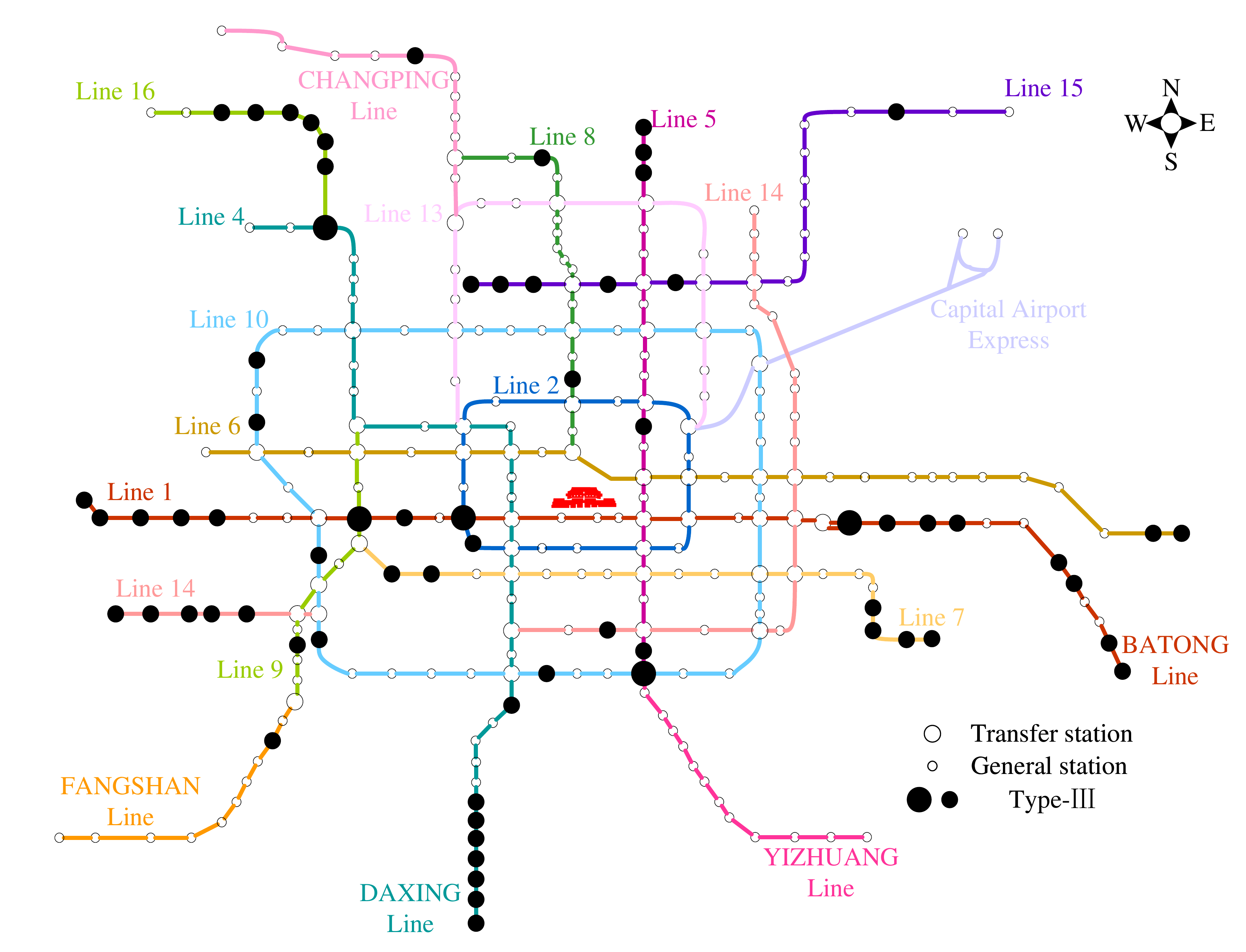}
}%
\subfigure[Type-\Rmnum{4}]{
\label{fig:f5.sub4}
\includegraphics[width=3.3in]{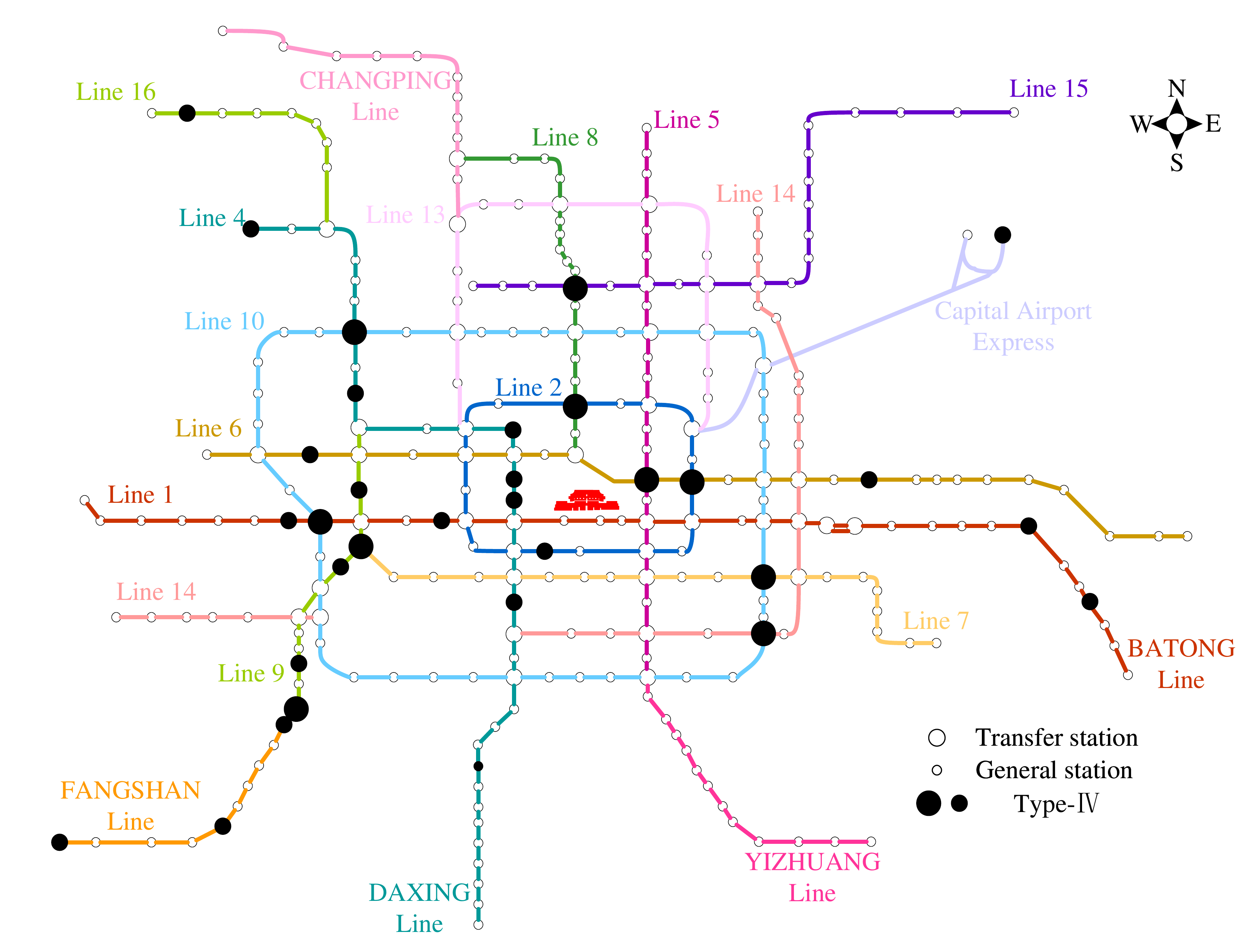}
}%

\centering
\caption{ A distribution map of different types of stations }
\label{fig:f5part}
\end{figure}

The four types of stations are marked in the Beijing URT line map, as shown in Figure \ref{fig:f5part}.
The spatial distribution characteristics and potential factors can be summarized from Figure \ref{fig:f5part} as follows.
\begin{itemize}
	\setlength{\itemsep}{0pt}
	\setlength{\parsep}{0pt}
	\setlength{\parskip}{0pt}

\item Type-\Rmnum{1}. The number of Type-\Rmnum{1}\ stations is the largest. In general, they spatially distribute on the circle line of the central area and connect branch lines in the outer areas of the city, and most of them are nontransfer stations represented by small solid black circles (see Figure \ref{fig:f5.sub1}). 
Those nontransfer stations are more affected by passenger flow than topology, and their importance changes significantly, especially in the morning and evening peaks.

\item Type-\Rmnum{2}. Stations have a major relationship with surrounding land types, many of which are near tourist attractions, airports and railway stations, as shown in Figure \ref{fig:f5.sub2}. 
Stations include both transfer and nontransfer stations, which have different topological characteristics such as degree and betweenness. 
However, every station has a high importance index at midday, which may be associated with the fact that such type of stations are usually with high travel demands at midday.

\item Type-\Rmnum{3}. Stations mainly distribute in the peripheral exurbs, many of which are concentrated at the end of branch lines. 
Therefore, passengers who have longer commuting distances from the suburbs to the city center usually start off from the station earlier.

\item Type-\Rmnum{4}. Stations consist mainly of transfer stations and surrounding stations, with a few distributed at the end of branch lines. 
They usually have a steady flow of passengers (staying high or low). 
Importance is influenced by both topology and passenger flow.

\end{itemize}

\begin{figure}[!ht]
    \centering
    \includegraphics[width=5in]{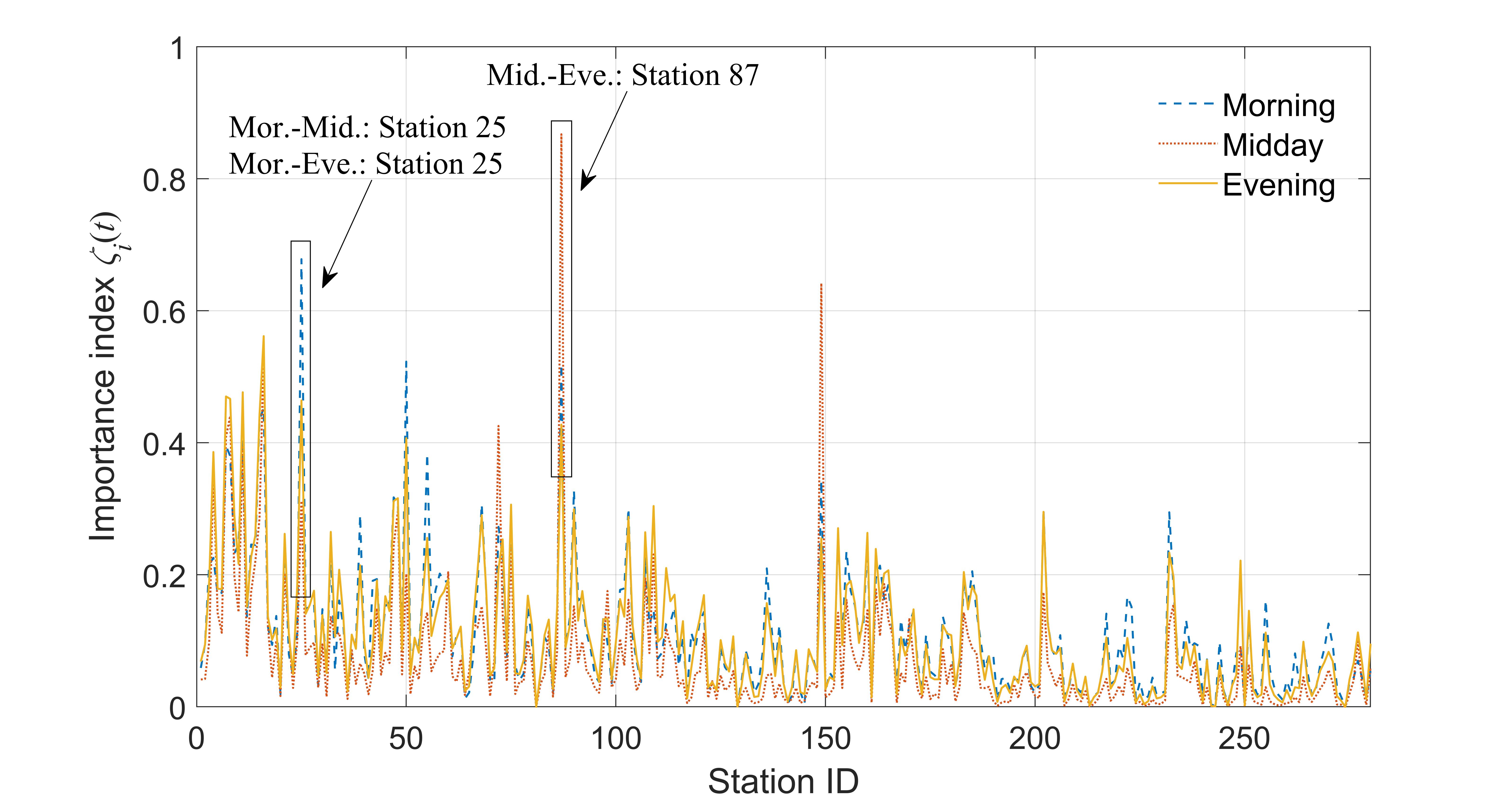}
    \caption{Importance index of each station at different time periods. The first 54 stations are transfer stations and then the nontransfer stations which are numbered in line order }
    \label{fig:f6}
\end{figure}

In these four types of curves, the most obvious peaks of importance indices appear mainly in the morning (see Figure \ref{fig:f4.sub1}), midday (see Figure \ref{fig:f4.sub2}) and evening (see Figure \ref{fig:f4.sub1} and \ref{fig:f4.sub3}). 
Therefore, we extract the typical periods in the morning, midday and evening of the day for further analysis. 
The importance index of each station in the morning-midday-evening periods can be obtained by using Equation \ref{eq4}. 
As shown in Figure \ref{fig:f6}, the importance index fluctuates from station to station over the same time period.
In particular, the difference in importance between adjacent stations is often large, reflecting that the importance of stations within the same line varies significantly. 
Therefore, even within the same line, different stations should be protected selectively according to the importance level. 
Moreover, most transfer stations (station IDs are less than 54) are with higher importance. 
A few nontransfer stations also have high indices of importance, since the importance of stations is not necessarily related to whether the station is a transfer station under the dual influence of passenger flow and topological structure.

\begin{figure}[!ht]
    \centering
    \includegraphics[width=6in]{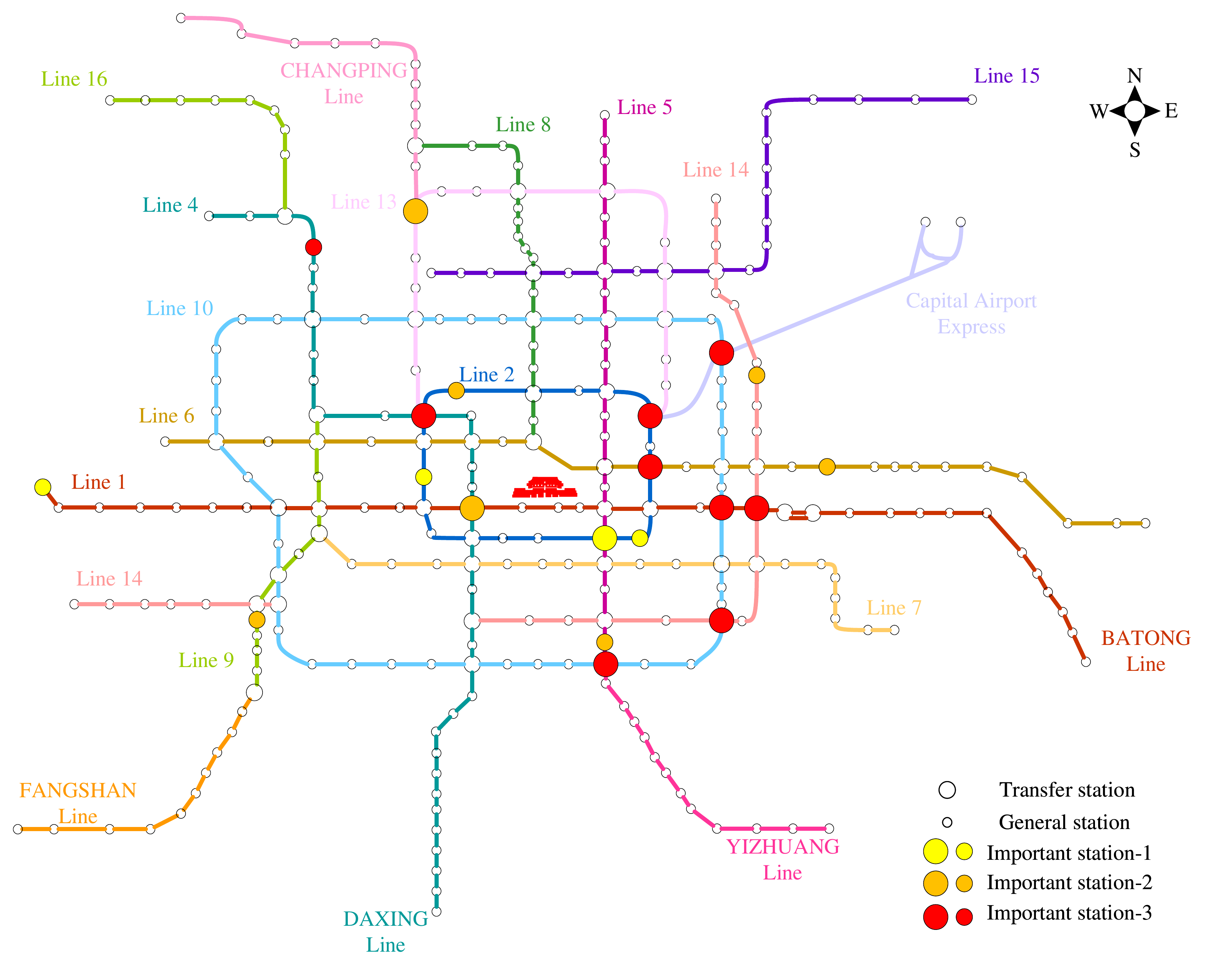}
    \caption{Distribution of high-importance stations at different levels. According to the frequency of one station appearing in the top 15 of the three periods (Table \ref{tab:my-table2}), the importance of those stations can be divided into three levels, i.e., the stations in Important Stations 1, 2 and 3 appear once, twice and three times in the morning-midday-evening periods, respectively. Obviously, the stations in Important Station-3 are relatively more important}
    \label{fig:f7}
\end{figure}

As shown in Figure \ref{fig:f6}, the index values of the same station also vary with time periods.
The maximum change in importance indices is found at Station 25 (difference values are 0.3662 and 0.2140) and Station 87 (difference value is 0.4396), caused mainly by the huge change in travel demand.
The top 15 stations ranked with respect to importance index $\zeta_i(t)$ are presented in Table \ref{tab:my-table2} for further analysis. 
It is found that most of the high-importance stations are transfer stations (represented by bold fonts), while the high-importance nontransfer stations are distributed mainly near the transfer stations or are the departure stations of major branch lines (see Figure \ref{fig:f7}). 
In addition, the high-importance stations in different time periods are not all the same but are more frequently repeated. 
For example, a total of 12 stations are all with high importance indices in the morning and evening periods. 
In conclusion, topological structure is still the main determinant of the importance level of URT stations, while passenger flow in the URT network is also found to be critical.

\begin{table}[!htbp]
\centering
\caption{The importance index rank of stations}
\scriptsize
\label{tab:my-table2}
\setlength\tabcolsep{6.6pt}
\begin{tabular}{lllllllllllll}
\toprule
\multirow{2}{*}{Rank} & \multicolumn{4}{c}{Morning} & \multicolumn{4}{c}{Midday}  & \multicolumn{4}{c}{Evening}  \\
 \cmidrule(r){2-5} \cmidrule(r){6-9} \cmidrule(r){10-13}
 & Station & $\zeta_i(t)$ & $B_{i}(t)$ & $C_{i}(t)$ & Station & $\zeta_i(t)$ & $B_{i}(t)$ & $C_{i}(t)$ & Station & $\zeta_i(t)$ & $B_{i}(t)$ & $C_{i}(t)$  \\
\midrule
1& \textbf{25} & 0.6780 & 0.0987 & 0.0065 & 87 & 0.8673 & 0.065 & 0.0273 & \textbf{16} & 0.5615 & 0.0627 & 0.0132 \\
2& \textbf{50} & 0.5229 & 0.0777 & 0.0229 & 149 & 0.6416 & 0.0999 & 0.0339 & \textbf{11} & 0.4764 & 0.144 & 0.0105 \\
3& 87 & 0.5192 & 0.0529 & 0.0172 & \textbf{16} & 0.5403 & 0.0682 & 0.0373 & \textbf{7} & 0.4701 & 0.1418 & 0.0105 \\
4& \textbf{16} & 0.4528 & 0.06 & 0.0219 & \textbf{8} & 0.4406 & 0.117 & 0.0424 & \textbf{8} & 0.4665 & 0.1179 & 0.0115 \\
5& \textbf{15} & 0.4345 & 0.1185 & 0.0612 & 72 & 0.4255 & 0.0406 & 0.0448 & \textbf{25} & 0.4641 & 0.0993 & 0.0454 \\
6& \textbf{11} & 0.4067 & 0.1438 & 0.0227 & \textbf{7} & 0.4084 & 0.1445 & 0.0410 & \textbf{15} & 0.4417 & 0.1258 & 0.0099 \\
7& \textbf{7} & 0.3953 & 0.1417 & 0.0247 & \textbf{11} & 0.3929 & 0.1538 & 0.0385 & 87 & 0.4276 & 0.0473 & 0.0220 \\
8& 55 & 0.3834 & 0.0143 & 0.0138 & \textbf{4} & 0.3556 & 0.1403 & 0.0562 & \textbf{50} & 0.4070 & 0.0709 & 0.0116 \\
9& \textbf{8} & 0.3800 & 0.1158 & 0.0282 & \textbf{25} & 0.3119 & 0.0942 & 0.0432 & \textbf{4} & 0.3860 & 0.1268 & 0.0138 \\
10& 149 & 0.3431 & 0.0884 & 0.0271 & \textbf{48} & 0.2955 & 0.0719 & 0.0520 & \textbf{48} & 0.3160 & 0.0755 & 0.0307 \\
11& 90 & 0.3277 & 0.0602 & 0.0117 & \textbf{15} & 0.2855 & 0.1266 & 0.0576 & \textbf{47} & 0.3116 & 0.0627 & 0.0167 \\
12& \textbf{47} & 0.3173 & 0.0639 & 0.0365 & \textbf{47} & 0.2447 & 0.0633 & 0.0651 & 75 & 0.3063 & 0.0485 & 0.0221 \\
13& 68 & 0.3068 & 0.0622 & 0.0634 & 75 & 0.2428 & 0.0511 & 0.0600 & 109 & 0.3042 & 0.0756 & 0.0317 \\
14& \textbf{48} & 0.3066 & 0.0801 & 0.0133 & 109 & 0.2308 & 0.0634 & 0.0727 & 90 & 0.2994 & 0.0637 & 0.1028 \\
15& 103 & 0.3007 & 0.0144 & 0.0182 & \textbf{14} & 0.2203 & 0.1276 & 0.0602 & 202 & 0.2953 & 0.0399 & 0.0184 \\
\bottomrule
\end{tabular}
\end{table}

\begin{table}[!htbp]
\centering
\caption{Importance ranking and related characteristic values for each operational line}
\scriptsize
\label{tab:my-table3}
\setlength\tabcolsep{26pt}
\begin{tabular}{llllll}
\toprule
\multirow{2}{*}{Time period} &\multirow{2}{*}{URT line} &\multirow{2}{*}{$\zeta_i(t)$} & \multicolumn{3}{c}{Characteristics} \\
 \cmidrule{4-6}
 &  &  & $B_{i}(t)$ & $C_{i}(t)$  & $D_{i}(t)$ \\
\midrule
\multirow{19}{*}{Morning} & Line 2 & 0.2201 & 0.0813 & 0.0609 & 265734 \\
 & Line 13 & 0.1996 & 0.0619 & 0.043 & 197822 \\
 & Line AE & 0.1926 & 0.0314 & 0.4549 & 94041 \\
 & Line 1 & 0.177 & 0.0808 & 0.1075 & 267821 \\
 & Line 5 & 0.1603 & 0.0808 & 0.0648 & 269721 \\
 & Line 10 & 0.141 & 0.0628 & 0.067 & 207198 \\
 & Line 4 & 0.1326 & 0.0674 & Inf & 224627 \\
 & Line 6 & 0.1323 & 0.0637 & 0.0603 & 210800 \\
 & Line FS & 0.1189 & 0.0389 & 0.0416 & 127107 \\
 & Line 9 & 0.1157 & 0.0889 & 0.0774 & 302880 \\
 & Line 14 & 0.0996 & 0.0406 & 0.0782 & 132499 \\
 & Line YZ & 0.0914 & 0.0316 & Inf & 101934 \\
 & Line CP & 0.0908 & 0.0211 & 0.1442 & 64345 \\
 & Line BT & 0.0863 & 0.0249 & 0.0811 & 78965 \\
 & Line 8 & 0.0812 & 0.037 & 0.1655 & 121284 \\
 & Line 7 & 0.0775 & 0.0388 & 0.1337 & 128208 \\
 & Line 15 & 0.0697 & 0.0301 & 0.1469 & 97785 \\
 & Line DX & 0.0459 & 0.0177 & 0.0829 & 57268 \\
 & Line 16 & 0.0305 & 0.0122 & 0.104 & 38257 \\
\midrule
\multirow{19}{*}{Midday} & Line AE & 0.1966 & 0.0343 & 0.3313 & 29134 \\
 & Line 2 & 0.1804 & 0.0846 & 0.105 & 81468 \\
 & Line 1 & 0.1344 & 0.0817 & 0.1198 & 79650 \\
 & Line 4 & 0.1277 & 0.0718 & Inf & 69976 \\
 & Line 13 & 0.1251 & 0.0513 & 0.1421 & 48314 \\
 & Line 5 & 0.0987 & 0.0794 & 0.1366 & 78350 \\
 & Line 9 & 0.0956 & 0.0899 & 0.1954 & 90028 \\
 & Line 6 & 0.0917 & 0.0623 & 0.1647 & 60788 \\
 & Line 10 & 0.0906 & 0.0586 & 0.1905 & 57150 \\
 & Line 14 & 0.0687 & 0.0374 & 0.3147 & 35796 \\
 & Line FS & 0.0575 & 0.0384 & 0.2176 & 36908 \\
 & Line YZ & 0.0551 & 0.0321 & Inf & 30616 \\
 & Line 7 & 0.0469 & 0.0375 & 0.404 & 36613 \\
 & Line BT & 0.0457 & 0.0231 & 0.344 & 21685 \\
 & Line CP & 0.04 & 0.0176 & 0.6999 & 16140 \\
 & Line 8 & 0.0312 & 0.0287 & 0.308 & 27926 \\
 & Line 15 & 0.0297 & 0.0243 & 0.4122 & 23350 \\
 & Line 16 & 0.0097 & 0.0093 & 0.6058 & 8608 \\
 & Line DX & 0.0084 & 0.0108 & 0.5841 & 10274 \\
 \midrule
\multirow{19}{*}{Evening} & Line 2 & 0.231 & 0.0816 & 0.0303 & 225076 \\
 & Line AE & 0.2184 & 0.0319 & 0.2265 & 79587 \\
 & Line 13 & 0.1959 & 0.0581 & 0.0562 & 156505 \\
 & Line 1 & 0.1859 & 0.0793 & 0.0447 & 221085 \\
 & Line 6 & 0.1564 & 0.0668 & 0.0849 & 186123 \\
 & Line 5 & 0.1521 & 0.0811 & 0.0556 & 228538 \\
 & Line 4 & 0.1456 & 0.0653 & Inf & 182911 \\
 & Line 10 & 0.1388 & 0.0599 & 0.0799 & 166892 \\
 & Line 14 & 0.1017 & 0.0395 & 0.1631 & 108629 \\
 & Line 9 & 0.1017 & 0.0829 & 0.0777 & 238347 \\
 & Line FS & 0.0927 & 0.0375 & 0.1592 & 103476 \\
 & Line YZ & 0.0915 & 0.0356 & Inf & 97289 \\
 & Line CP & 0.0861 & 0.0207 & 0.6104 & 53359 \\
 & Line 7 & 0.0815 & 0.0378 & 0.1888 & 105210 \\
 & Line BT & 0.0782 & 0.0235 & 0.2372 & 62592 \\
 & Line 8 & 0.0743 & 0.0345 & 0.1269 & 95371 \\
 & Line 15 & 0.0571 & 0.0277 & 0.1376 & 76006 \\
 & Line DX & 0.0279 & 0.0154 & 0.4598 & 42121 \\
 & Line 16 & 0.0221 & 0.0102 & 0.4044 & 26957 \\
 \bottomrule
\end{tabular}
\begin{tablenotes}\scriptsize
\item[1] Note: In the actual data, we found that Xinjiekou Station of Line 4 and Xiaocun Station of Yizhuang Line had no card swiping record on that day. Therefore, the demand closeness centrality $C_{i}(t)$ of these two lines is infinite.
\end{tablenotes}
\end{table}

Based on the importance index of each station, the importance ranking of each operational line can be obtained, which helps to understand the average importance of each line in the whole Beijing URT network. 
As shown in Table \ref{tab:my-table3}, the importance of lines at different time periods is ranked by average importance index $\zeta_i(t)$). 
The results indicate that Line 2 is with the highest importance index in the morning and evening, while the Airport Line is the highest at midday, consistent with previous studies that airports are important nodes in a URT network \citep{Wang2020}.
Some other important lines include branch Line 13 and diameter Line 1. 
The traditional topological characteristics of these lines, such as flow betweenness, demand closeness centrality and intensity, are also high, reflecting the effectiveness of our proposed importance index. 
Some of the peripheral lines, such as the Yizhuang Line, are also of high importance, especially in the morning and evening periods. 
In addition, we observe a large difference in importance between the two ring lines, i.e., Lines 2 and 10. 
Although the passenger flow of Line 10 is larger than the passenger flow of Line 2 (as shown in Figure \ref{fig:f3}), the importance of Line 2, which is closer to the center, is higher, showing that the importance of ring lines is closely related to their location in the city and topology in a URT network. 
As a reference metric, the flow betweenness and other common indices can only represent limited features of station importance, and thus the corresponding ranking results are different from that resulted from the proposed importance index.

\subsection{Network vulnerability under short-delay disruptions}

We investigate the vulnerability of the Beijing URT network under short-delay disruptions, i.e., the delay $\tau_\text{disr} \leq \tau_\text{disr}^*$ where $\tau_\text{disr}^*=60$ minutes considering the statistics of travel time and passengers' tolerance of delay.
The delays of 5, 10, 20, 40 and 60 minutes are employed to incorporate the different impacts of operation incidents on passenger demands and network topology. 
The top 15 high-importance stations identified in the morning period (see Table \ref{tab:my-table2}) are selected to be disturbed to imitate short-delay disruptions.
The vulnerability metric $\psi_\text{short} (\textbf{G}, t,\tau_\text{disr})$ is calculated by using Equation \ref{eq5}, which represents the average delayed passenger demands. 
The higher the value $\psi_\text{short} (\textbf{G}, t,\tau_\text{disr})$ is, the more vulnerable the network is.

\begin{figure}[!ht]
\centering
\subfigure[Changes in performance at different stations]{
\includegraphics[width=3in]{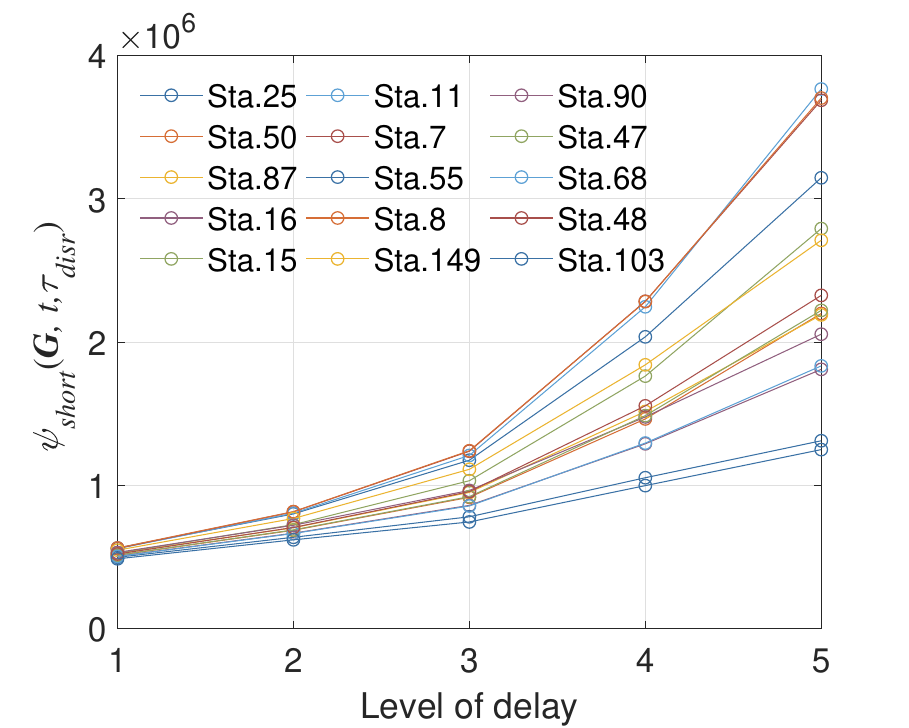}
\label{fig:f8.sub.1}
}
\subfigure[ Changes in performance at different delays]{
\includegraphics[width=3in]{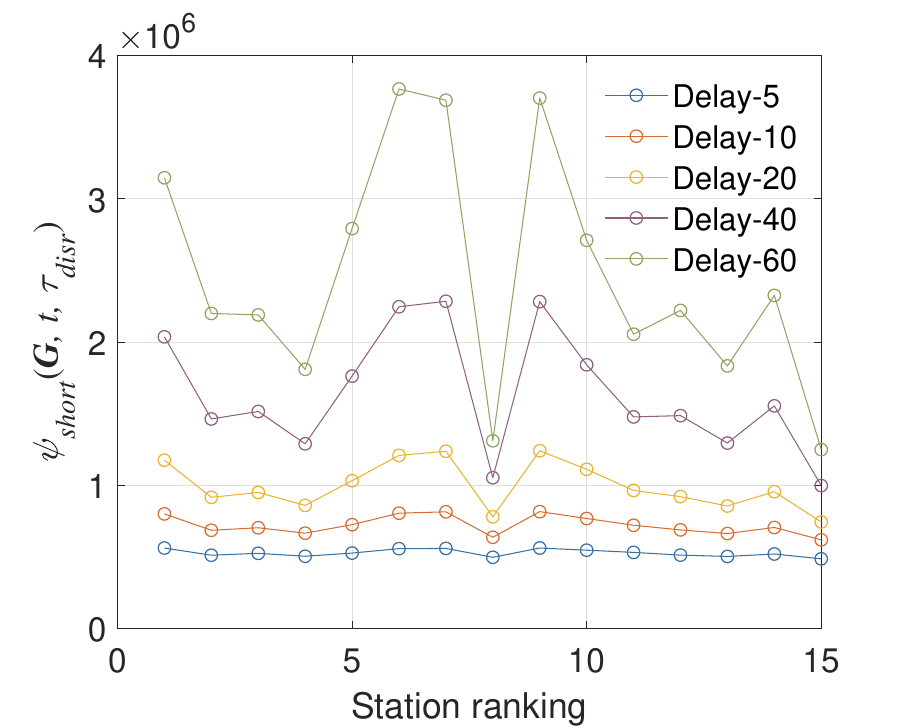}
\label{fig:f8.sub.2}
}
\caption{Network vulnerability changes under different delays during the morning period}
\label{fig:f8}
\end{figure}

It can be seen from the experimental results presented in Figure \ref{fig:f8} that the impact of stations on network vulnerability also changes with the increase of delay. 
However, the station importance is not linearly related to the degree of vulnerability. 
The stations with high-importance, such as the 2nd, 3rd and 4th stations, may have a diminishing impact on network vulnerability as delay increases. 
Some of the lower-ranked stations can have a greater impact as delays increase (such as the 9th and 14th stations), which may be related to the increase in travel demand affected by operational delay. 
As the delay increases, the affected travel demands vary for different stations, and the decline of network performance shows a nonlinear change. 
The experiments at midday and evening periods also have similar results.

\begin{figure}[htbp]
\centering
\subfigure[Cluster]{
\includegraphics[width=2.8in]{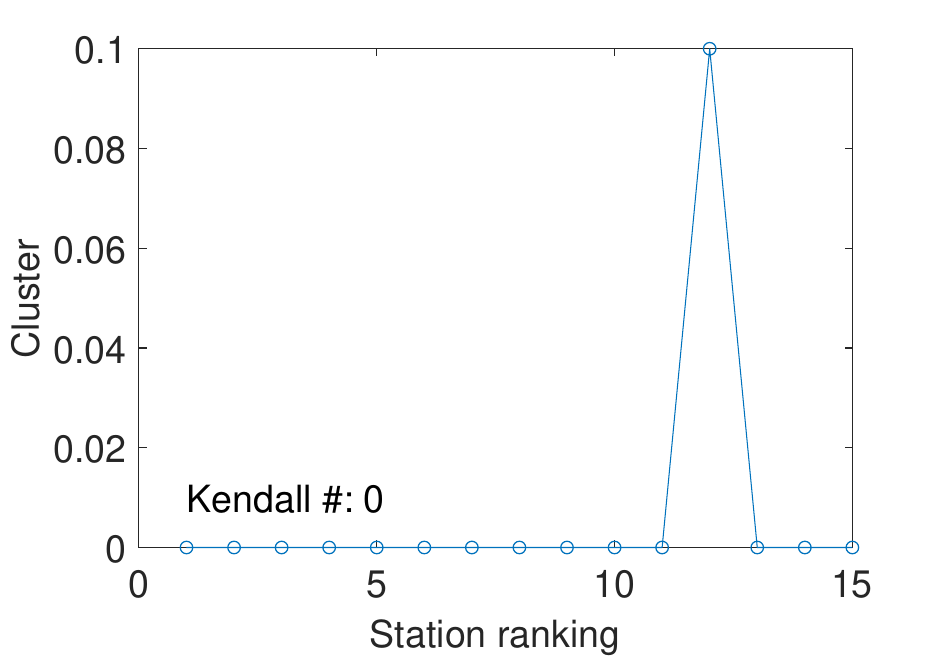}
\vspace{-4 mm}
}%
\subfigure[Degree]{
\includegraphics[width=2.8in]{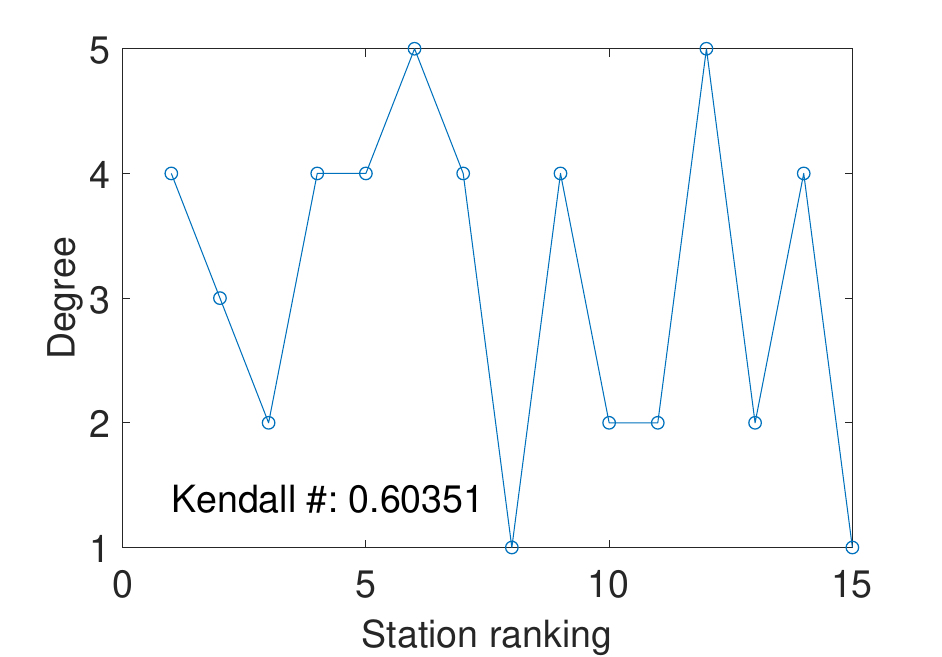}
\vspace{-4 mm}
}%
\vspace{-4 mm}

\subfigure[Betweenness]{
\includegraphics[width=2.8in]{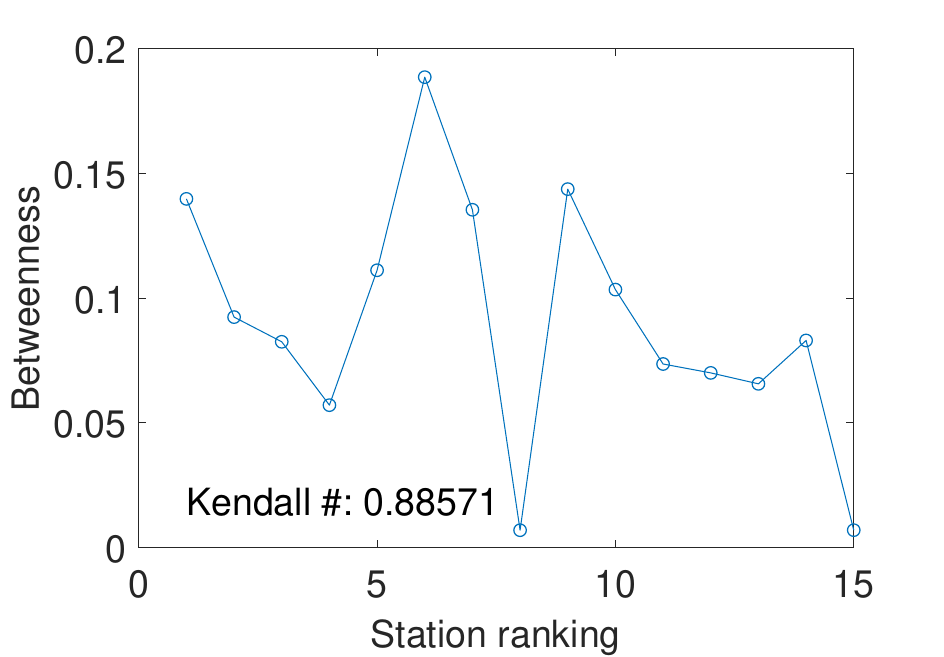}
\vspace{-4 mm}
}%
\subfigure[Closeness]{
\includegraphics[width=2.8in]{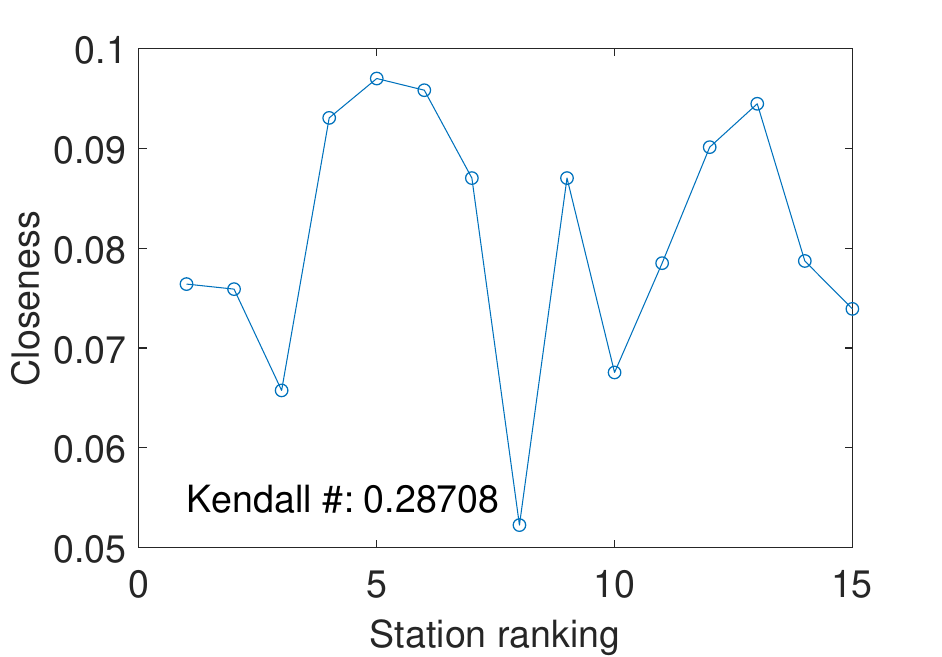}
\vspace{-4 mm}
}%
\vspace{-4 mm}

\subfigure[Number of passengers entering]{
\includegraphics[width=2.8in]{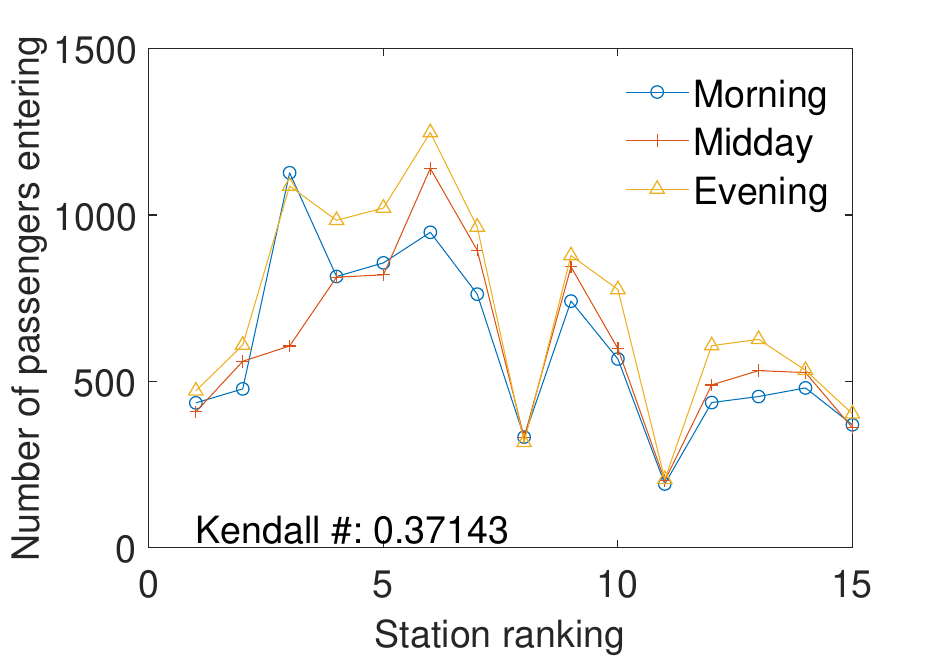}
\vspace{-4 mm}
}%
\subfigure[Intensity]{
\includegraphics[width=2.8in]{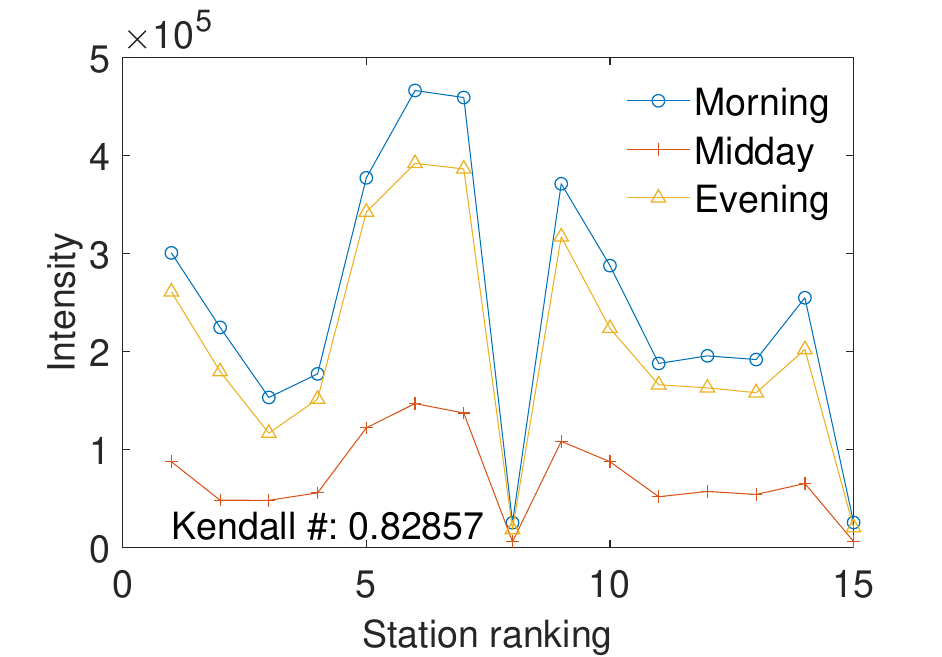}
\vspace{-4 mm}
}%
\vspace{-4 mm}

\subfigure[Flow betweenness]{
\includegraphics[width=2.8in]{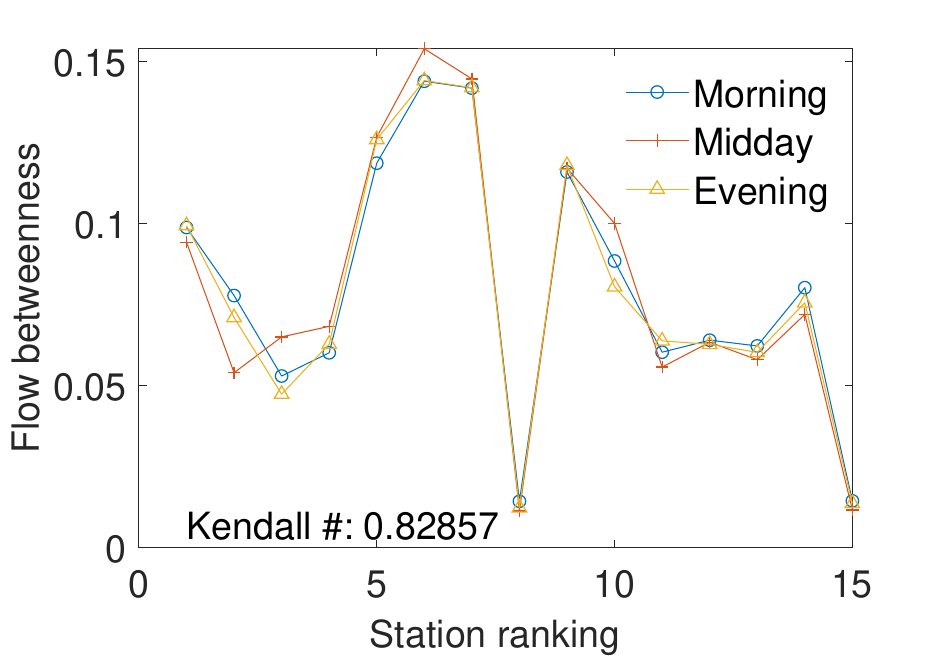}
\vspace{-4mm}
}%
\subfigure[Demand closeness centrality]{
\includegraphics[width=2.8in]{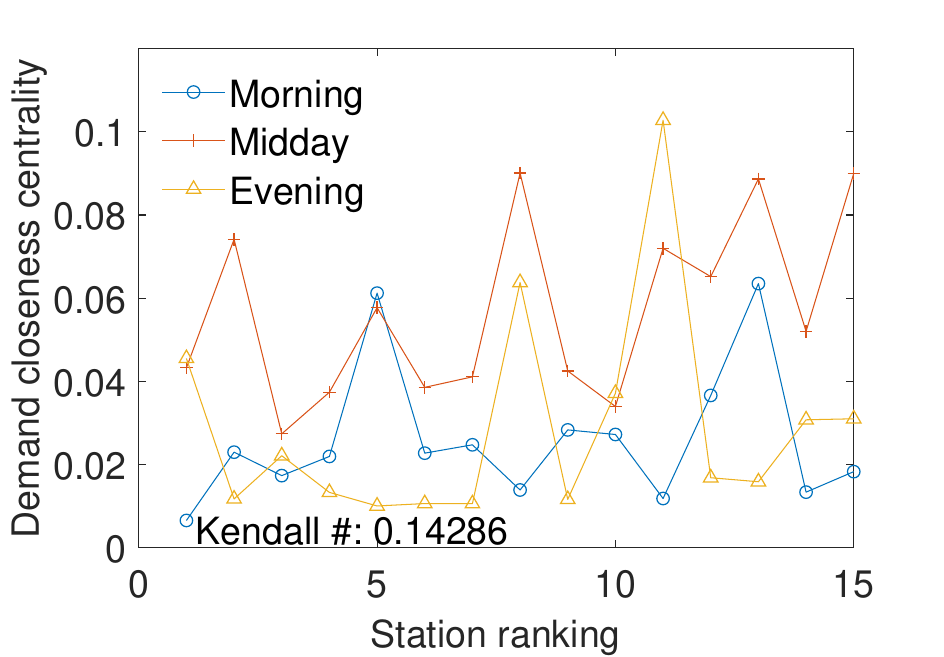}
\vspace{-4mm}
}%
\vspace{-3mm}
\centering
\caption{Eight typical characteristics of the station. Kendall: the correlation coefficient between characteristic curves (when there are different time periods, the curve of the morning time period is adopted) and the vulnerability curve }
\label{fig:f9}
\end{figure}

To further explore the reasons why stations have different impacts on vulnerability over time, we analyze the topology and passenger flow characteristics of those stations in Figure \ref{fig:f9}. 
We conduct the Kendall correlation analysis between the vulnerability curve (see Figure \ref{fig:f8.sub.2}) of 60-minute delay and the characteristic curves of the stations. 
The topological betweenness is shown to have the strongest correlation with the vulnerability curve.
The intensity and flow betweenness of stations that are related to passenger flow also have a high correlation, indicating that some stations that serve as bridges (i.e., the reasonable paths that link the origin station and destination station) and have a large number of passengers usually have the greatest influence on network vulnerability.

The correlation between vulnerability curves and the selected four characteristics curves with high correlation under different time periods and delays is further analyzed in Table \ref{tab:delay}.
In the same period, the correlation between curves increases gradually with the increasing delay.
With a longer delay, the characteristics of the station itself influence the network vulnerability more apparently, which can lead to the changing trend of network performance.
Notably, remarkable characteristics are not all the same in different time periods, as shown in bold fonts in Table \ref{tab:delay}.
Multiple characteristics are known to jointly determine the change in network vulnerability.
However, the key factors affecting network vulnerability may be different under different degrees of disturbance or stages of disruption.

\begin{table}[htbp]
\centering
\caption{Correlation analysis under different delay degrees at different time periods}
\scriptsize
\setlength\tabcolsep{18.5pt}
\begin{tabular}{llllll}
\toprule
\multirow{2}{*}{Time period} & Delays  & \multicolumn{2}{c}{Topological characteristics} & Intensity   & Flow betweenness  \\
\cmidrule{3-4}
 & (min) &  Degree &  Betweenness & $D_{i}(t)$ & $B_{i}(t)$  \\
\midrule
\multirow{3}{*}{Morning} & 5 & 0.38797 & 0.73333 & 0.60000 & 0.60000 \\
 & 20 & 0.45263 & 0.77143 & 0.67619 & 0.67619 \\
 & 60 & 0.60351 & \textbf{0.88571} & 0.82857 & 0.82857 \vspace{1mm} \\
\multirow{3}{*}{Midday} & 5 & 0.33635 & 0.71429 & 0.58095 & 0.61905 \\
 & 20 & 0.31232 & 0.71429 & 0.69524 & 0.65714 \\
 & 60 & 0.45647 & 0.79048 & \textbf{0.80952} & 0.73333 \vspace{1mm}\\
\multirow{3}{*}{Evening} & 5 & 0.39977 & 0.65714 & 0.48571 & 0.46667 \\
 & 20 & 0.35408 & 0.65714 & 0.67619 & 0.65714 \\
 & 60 & 0.51399 & 0.69524 & \textbf{0.75238} & 0.73333 \\
 \bottomrule
 \label{tab:delay}
\end{tabular}
\end{table}

\subsection{Network vulnerability under long-delay disruptions}

Assume that the station is completely shut down and all transportation functions are disabled when $\tau_\text{disr} > 60$ minutes, and then single station failure and interval failure are simulated. 
The operational efficiency $\psi_\text{long} (\textbf{G}', t,\tau_\text{disr})$ in Equation \ref{eq6} is selected to be the vulnerability metric.

\begin{figure}[!ht]
    \centering
    \includegraphics[width=3.2in]{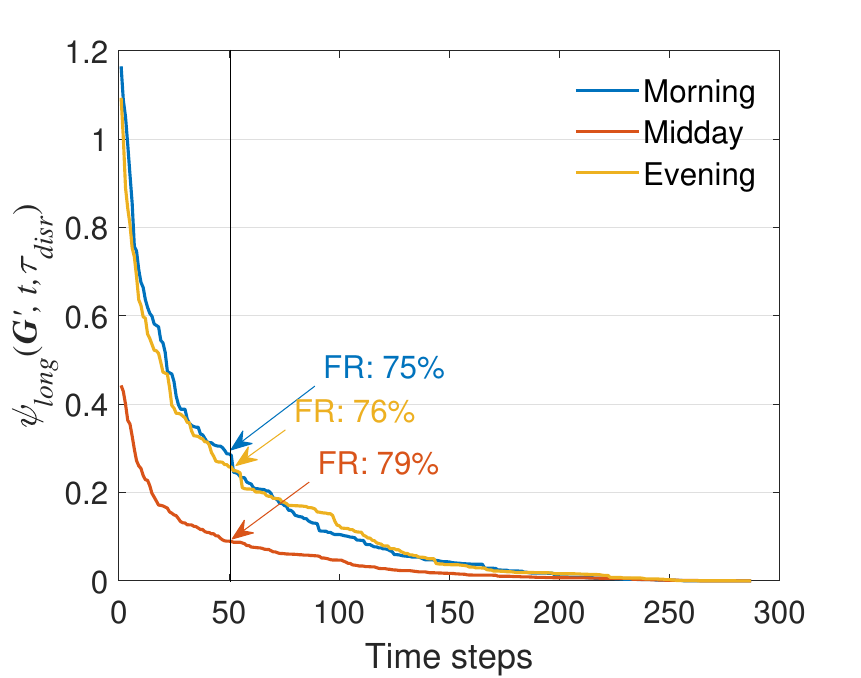}
    \caption{The changes of network vulnerability under station failure. FR: failure rate of stations in the Beijing URT network}
    \label{fig:f10}
\end{figure}

First, we investigate the vulnerability performance variation of the Beijing URT network under single station failure.
Based on the importance ranking of stations identified in Section \ref{section4.2}, the network vulnerability is simulated by means of deliberate attack. 
The results are shown in Figure \ref{fig:f10}. 
For the morning-midday-evening periods, the vulnerability index of the Beijing URT network caused by station failure shows an exponential decline with an obvious heavy tail distribution.
Compared with the morning and evening periods, the network vulnerability at midday is at a lower initial index, while it decreases the most by 79\% when the time steps reach 50. 

\begin{figure}[!ht]
    \centering
    \includegraphics[width=3.2in]{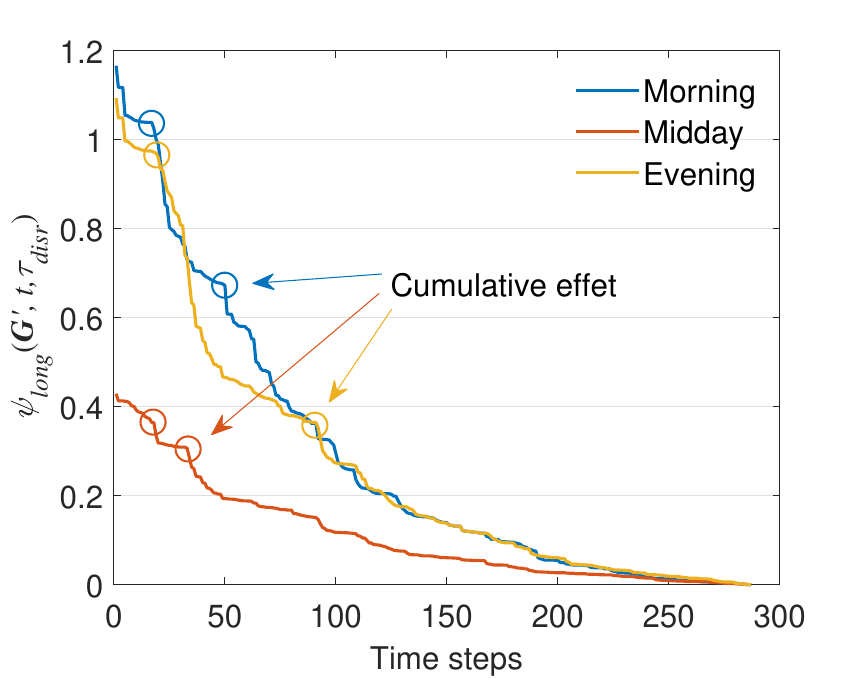}
    \caption{The changes of network vulnerability under interval failure}
    \label{fig:f11}
\end{figure}

Second, we simulate the continuous failure of multiple stations from the highest-importance operational line to the lowest-importance one in Table \ref{tab:my-table3}.
Within a line, let the highest-importance station fail first and then let the failure propagate along the direction of the second highest-importance station.
The simulation results are shown in Figure \ref{fig:f11}. 
Compared with single station failure, the continuous failure has similar time steps to completely crash the network, while the detailed process of performance variation is different.
In the single station failure, the URT network is close to breakdown when time steps reach 200, which is earlier than continuous failure when the time steps need to reach 250. 
There is also a special cumulative effect.
The decline of network performance is obviously segmented and stratified, and network performance will fall off a cliff when the cumulative number of failed stations reaches a boundary point.
This is probably because the failure of multiple stations in a single line has less effect on network vulnerability. 
The failure of stations connecting different lines, especially the boundary stations of each layer in the vulnerability curve, is the main cause of the performance decline, and some important stations within the line will act as catalysts to accelerate the performance decline of the Beijing URT network.

\begin{figure}[!ht]
\centering
\subfigure[Direction towards lower importance station]{
\includegraphics[width=3in]{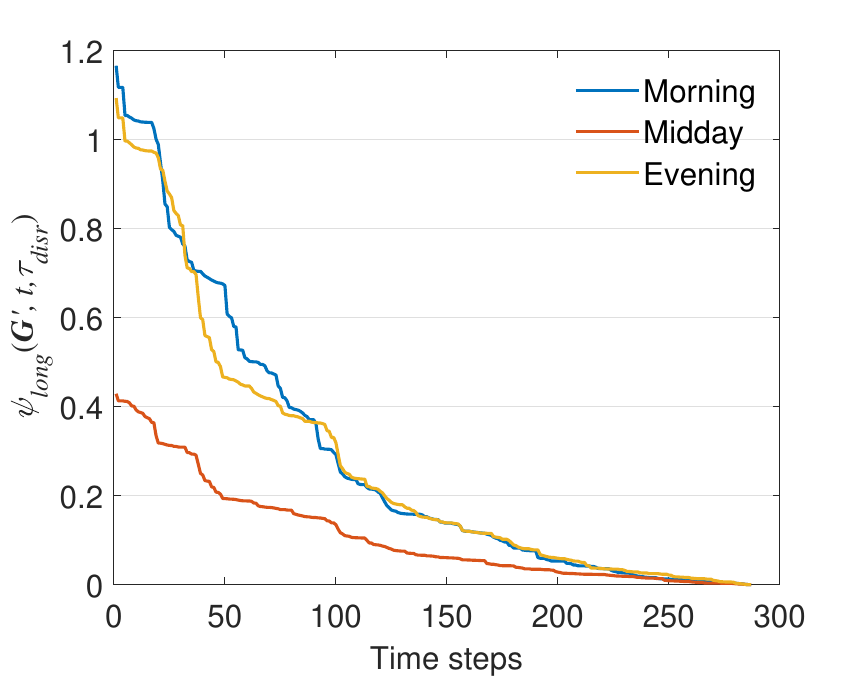}}
\subfigure[Both directions]{
\includegraphics[width=3in]{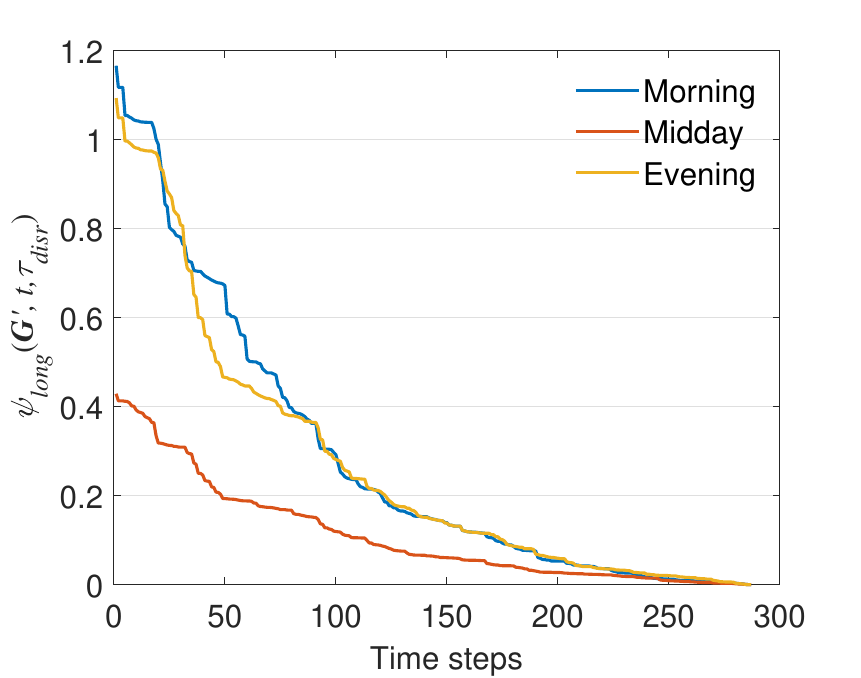}}
\caption{Interval failure in different directions}
\label{fig:f12}
\end{figure}

To explore the influence of the failure direction within a line on the URT network vulnerability, we carry out experiments in opposite directions of failure and two directions of failure, respective. 
The results in Figure \ref{fig:f12} show that the failure direction has little influence on the vulnerability. 
Instead, the decline of network performance depends mainly on the failure of critical lines or stations. 
Combined with Figures \ref{fig:f11} and \ref{fig:f12}, the significant decline in vulnerability index $\psi_\text{long} (\textbf{G}', t,\tau_\text{disr})$ occurs mainly after the failure of several specific stations or the cumulative failure of multiple stations. 
In addition, high ranking of importance of a single line or station, according to the Section \ref{section4.2}, does not mean that it has a great influence on the network vulnerability. 
The degree of influence on the network vulnerability is often the cumulative effect of multiple stations or lines. 
Therefore, it is necessary to avoid the continuous failure of multiple stations, especially the failure of an entire line and connected lines.

\begin{figure}[!ht]
\centering
\subfigure[Simultaneous failure of two stations]{
\includegraphics[width=3in]{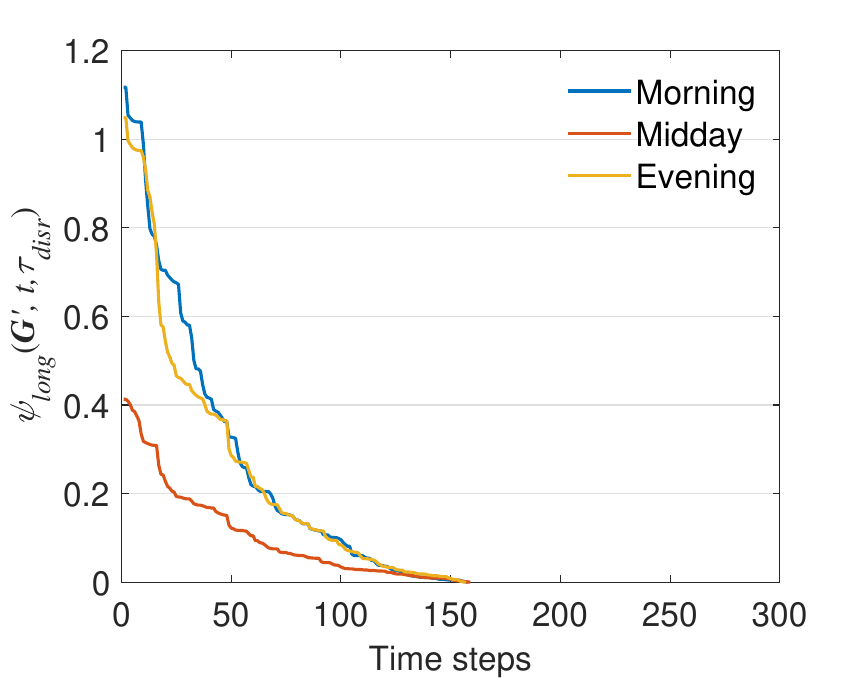}
}
\subfigure[Simultaneous failure of three stations]{
\includegraphics[width=3in]{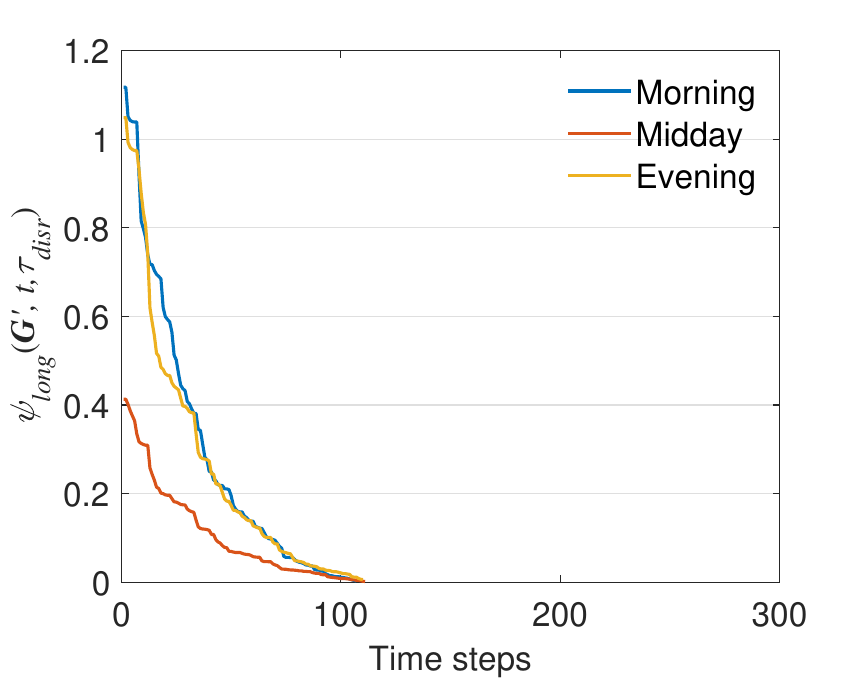}
}
\centering
\caption{Interval failure within the same line }
\label{fig:f13}
\end{figure}

At last, we focus on the influence of special interval failure (i.e., two or three adjacent stations fail simultaneously) on URT network vulnerability.
The simulation method of interval failures is similar to that of continuous failure of multiple stations.
The simulation starts from the highest-importance operational line to the lowest-importance one in Table \ref{tab:my-table3}, and within a line the highest-importance interval fails first. 
The results are similar as shown in Figure \ref{fig:f13}. 
Interval failures within the same line degrade network performance faster than station failures.
The segmented and stratified decline of the vulnerability curve is still evident, indicating that the dominant loss of the URT network performance is due to the failure of certain critical stations and lines.

\begin{figure}[htbp]
    \centering
    \includegraphics[width=3.2in]{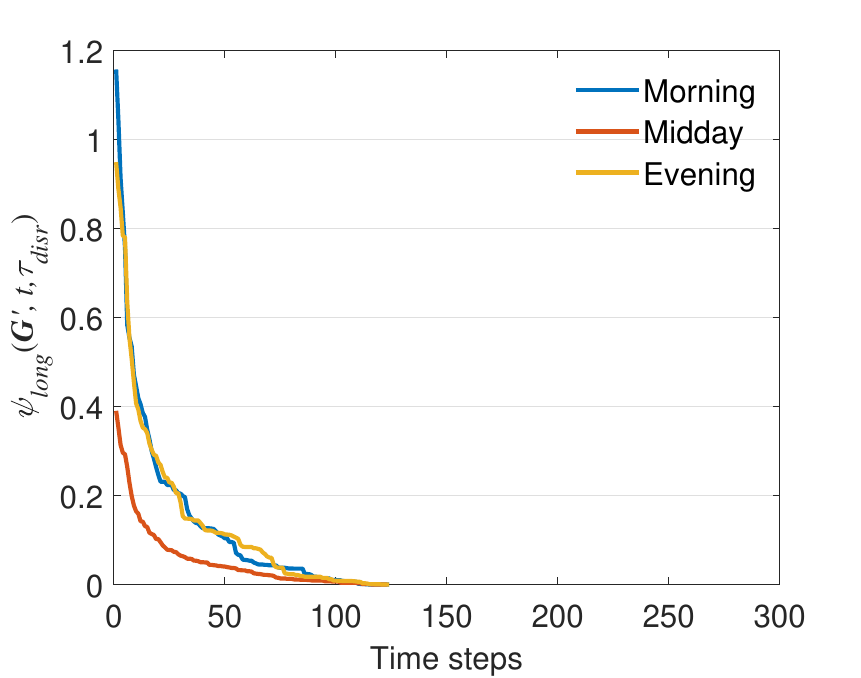}
    \caption{Interval failure among different lines}
    \label{fig:f14}
\end{figure}

Then, we conduct interval failure experiments where intervals fail randomly among lines. 
According to the ranking of the total importance of any two adjacent stations, the interval failure composed of two stations is tested. 
The results in Figure \ref{fig:f14} show that the network performance degrades faster in the case of interval failure of different lines so that, in practice, interval failures between different lines should be avoided as much as possible. 
Confining the consequences of disturbance to a single interval or line would greatly avoid a rapid breakdown of the entire network and race against time for rescue operations.

\section{Discussion and conclusion}\label{section5}

Considering the time-varying travel demands, this paper proposes an accessibility-based evaluation method to analyze the dynamic vulnerability of a URT system from the network perspective. 
Based on the static network topology and dynamic passenger flow, the importance index of stations and the vulnerability metric are established. 
By taking the URT system in Beijing, China as an example, this paper analyzes the impact of different disruption types disturbed by daily operation incidents on network vulnerability. 
New findings are summarized as follows.
\begin{itemize}
	\setlength{\itemsep}{0pt}
	\setlength{\parsep}{0pt}
	\setlength{\parskip}{0pt}	
	
\item The distribution of high-importance stations varies with time periods as there are many recurring high-importance stations, showing that the travel demand as the one who uses the URT system is an important influencing factor.
The high-importance stations include not only the transfer stations with high passenger flow but also the nontransfer stations around them and the departure station.

\item Moreover, it is found that some high-importance stations exist independently of the dynamic passenger demands, indicating that the topological structure plays a decisive role. 
Comparing with the high-importance stations, the ranking of high-importance lines is more stable over different time periods.

\item 
In the case of short-delay disruption, the impact of high-importance stations on network vulnerability varies with the increase in delayed passenger demand. 
The stations that have the greatest impact on network vulnerability are those that act as bridges and have large passenger flows, while not all these stations are hubs.
In the case of long-delay disruption, the vulnerability of Beijing URT decreases exponentially and shows an obvious heavy tail distribution.

\item The degradation of network performance under continuous interval failure is obviously segmented and stratified. 
The interval failure within the same line has little effect on the vulnerability of Beijing URT network, while the failure between different lines and their cumulative effect is the main reason for the sharp decrease in network performance. 
Some high-importance stations or intervals within the line can act as catalysts to accelerate the performance degradation of the URT network, regardless of the direction of the station failure.

\end{itemize}

In practice, the resources for protection and rescue are limited, and URT managers should selectively protect and repair important components (stations, intervals and lines), which will greatly reduce system breakdown probability and recovery time. 
In preventing disruptions, a reasonable scheme of resource allocation should be carried out for different components in different time periods. 
Some recurring components of importance identified in different periods need us to pay special attention, and the adjacent nontransfer stations should also be protected. 
In addition, components for which indices of importance vary greatly in different periods tend to be more critical than the components constantly with high-importance. 
Once a disruption occurs, the rescue force should pay more attention to the key stations at the boundary of URT lines or intervals, i.e., the end nodes of each layer in a step-down vulnerability curve. 
If the key stations at the boundary of URT lins or intervals fail, the network performance will fall off a cliff. 
In the incident of continuous station failures, multiple interval failures between different lines should be avoided as much as possible, because the multiple interval failures might cause the network performance to plummet and crash the entire URT system.

The assumption that passengers will wait when the delay time increases is not exactly consistent with reality. 
In the future, field surveys or questionnaires might be consider to obtain more accurate travel data for the vulnerability assessments.
In the clustering process of the importance curves, we find that the results are influenced by time granularity, which is expected to improve. 
It is also expect to see more discussion regarding the influence of more types of interval failures on the vulnerability of a URT network, including the division of intervals within a line. 
The distribution of high-importance stations is also related to the geographical environment, population density and economic activities, and the correlation can be further explored by combining relevant data such as points-of-interest data . 
In addition, how to effectively recover the URT system performance after disruptions is also an important research direction in the future.

\section*{Acknowledgment}

The research is funded by 
National Natural Science Foundation of China (71871010).

\bibliographystyle{model2-names}
\bibliography{library}

\end{spacing}
\end{document}